\documentclass[preprint,review,12pt]{elsarticle}

\usepackage{amsthm}
\usepackage{hyperref}
\usepackage{moreverb}
\usepackage{graphicx,subfig,wrapfig}
\usepackage{algorithm}
\usepackage{algpseudocode}
\usepackage{color}
\usepackage{amsmath}
\usepackage{amssymb}
\usepackage{bm}
\usepackage{afterpage}
\usepackage[margin=1in]{geometry}
\usepackage{mathrsfs}
\usepackage{tikz}
\usepackage{pstricks}
\usepackage{diagbox}
\usepackage{arydshln}
\usepackage{multirow}
\usepackage{appendix}

\newcommand{\bx}{\boldsymbol{x}}

\newcommand{\btheta}{\mathbf{\Theta}}
\newcommand{\pL}{\mathcal{L}}
\newcommand{\pB}{\mathcal{B}}
\newcommand{\mS}{\mathrm{S}}

\newcommand{\mR}{\mathbb{R}}
\newcommand{\mN}{\mathbb{N}}
\newcommand{\bg}{\boldsymbol{\gamma}}
\newcommand{\bB}{\boldsymbol{B}}
\newcommand{\pOmega}{\partial \Omega}

\newcommand{\pjoi}{\partial_{j} \Omega_{i}}
\newcommand{\pexoi}{\partial_{ex}\Omega_{i}}

\def\linteri{\Lambda_{i}}

\def\lr{\color{black}}

\usepackage{algorithm,algpseudocode}
\usepackage{amsmath}
\usepackage{txfonts}
\usepackage{afterpage}

\usepackage{lipsum}
\usepackage{amsfonts}
\usepackage{graphicx}
\usepackage{epstopdf}

\usepackage{amsopn}

\journal{}

\begin{document}

\begin{frontmatter}

\title{A deep domain decomposition method based on Fourier features} 
\tnotetext[t1]{This work is supported by the National Natural Science Foundation of China (No. 12071291), the Science and Technology Commission of Shanghai Municipality (No. 20JC1414300) and the Natural Science Foundation of Shanghai (No. 20ZR1436200).
}
\author{Sen Li}
\ead{lisen@shanghaitech.edu.cn}
\author{Yingzhi Xia}
\ead{xiayzh@shanghaitech.edu.cn}
\author{Yu Liu}
\ead{liuyu@shanghaitech.edu.cn}
\author{Qifeng Liao\corref{cor1}}
\ead{liaoqf@shanghaitech.edu.cn}
\cortext[cor1]{Corresponding author}
\address{School of Information Science and Technology, ShanghaiTech University, Shanghai 201210, China }

\begin{abstract}
In this paper we present a Fourier feature based deep domain decomposition method (F-D3M) for 
partial differential equations (PDEs). Currently, deep neural network based methods are actively developed for solving PDEs, but their efficiency can  degenerate for problems with high frequency modes. 
In this new F-D3M strategy, overlapping domain decomposition is conducted for the spatial domain, such that
high frequency modes can be reduced to relatively low  frequency ones. 
In each local subdomain, multi Fourier feature networks (MFFNets) are constructed,
where efficient boundary and interface treatments are applied for the corresponding loss functions.
We present a general mathematical framework of F-D3M, validate its accuracy and 
demonstrate its efficiency with numerical experiments. 

\end{abstract}


\end{frontmatter}


\section{Introduction}\label{sec:intro}
During the past few years, there has been a rapid development in deep learning based methods for 
solving partial differential equations (PDEs). 
This explosion in interest has been driven by the need of surrogate modeling for uncertainty quantification, 
efficient treatments for 
complex geometries, and meshless approximation for
high-dimensional problems. 
The main idea of these methods is to reformulate a PDE problem as an optimization problem 
and train deep neural networks to approximate the solution through minimizing the corresponding loss function. 
These deep learning based methods typically include, 
deep Ritz methods \cite{weinan2017proposal,weinan2018deep}, 
physical informed neural networks \cite{raissi2019physics, karniadakis2021physics}, 
 deep Galerkin methods \cite{sirignano2018dgm}, 
Bayesian deep convolutional encoder-decoder networks \cite{zhu2018bayesian,zhu2019physics}
and 
deep multiscale model learning \cite{wang2020deep}. 
In addition, 
deep neural network methods for complex geometries and interface problems are proposed in 
\cite{sheng2021pfnn, gao2021phygeonet, wang2020mesh}, and 
deep mixed residual methods are developed for high-order PDEs in 
\cite{lyu2022mim}.

Whilst achieving great successes, there still exist challenges for these new methods,
especially when the underlaying problem has high frequency modes. 
As addressed in   \cite{xu2020frequency,luo2019theory,rahaman2019spectral,ronen2019convergence},
deep neural networks have the pathology properties of the frequency principle and the spectral bias, 
which show that neural networks are typically efficient  for fitting objective functions with low frequency modes
but can be inefficient for high frequency functions. 
To improve the efficiency for solving problems with high frequencies,
multi-scale networks are proposed through radial scaling in frequency domains \cite{cai2019multi,li2020multi,wang2020multi,liu2020multi},
where learning a high frequency function is converted to  learning a lower frequency one by 
scaling the input variables with multiple appropriate factors. 
In  \cite{tancik2020fourier}, it is shown that through mapping input coordinates with a simple Fourier feature 
can alleviate the spectral bias of fully connected neural networks,
where the  Neural Tangent Kernel theory \cite{jacot2018neural} is used. 
Moreover, effective multi-scale Fourier feature networks are proposed based on the random Fourier features
in \cite{wang2021eigenvector}. 

 To further improve the computational efficiency of neural network based approaches,
 this paper is devoted to domain decomposition methods. 
Due to their scalability from the divide and conquer principle, 
 domain decomposition methods are widely used for solving complex PDE problems. 
The books  by Chan and Mathew \cite{chan1994domain}, Quarteroni and Valli \cite{quarteroni1999domain}, and Toselli and Widlund \cite{toselli2004domain}, give a general review. 
Domain decomposition for neural network based methods currently gains a lot of interests,
and new methods are actively developed. 
In \cite{li2019d3m}, we propose a deep domain decomposition method 
(D3M), where neural networks are constructed in subdomains with the deep Ritz method. Li et al. \cite{li2020deep} propose physics-informed neural networks based on domain decomposition  for elliptic problems. 
In \cite{jagtap2020conservative}, Jagtap et al. develop conservative physics-informed neural networks for 
conservation laws.  
Through decomposing test function spaces, Kharazmi et al.\ \cite{kharazmi2021hp}  develop variational physics-informed neural networks through 
projection onto spaces of high-order polynomials. 
In  \cite{dong2021local}, Dong and Li introduce local extreme learning
machines with local field solutions represented by feed-forward neural networks. 
Besides, Moseley et al.\ \cite{moseley2021finite} propose 
finite basis physics-informed neural networks with separate input normalization  over subdomains, and Sheng and Yang \cite{sheng2022pfnn2} develop a penalty-free neural network method based on domain decomposition.
 
In this paper we present a new Fourier feature based deep domain decomposition method
(F-D3M) to curb the difficulty posed by high frequency modes. The method employs overlapping domain 
decomposition such that high frequency problems can be efficiently solved locally in subdomains. 
A distinct feature of F-D3M is that relative frequencies are reduced along with decomposing 
a global domain  into small subdomains.
In each local subdomain, new multi Fourier feature networks (MFFNets) are introduced with loss functions
focusing on interior parts of local solutions only. 
While boundary conditions for both exterior boundaries and 
domain decomposition interfaces are not involved in these loss functions, each MFFNet can be efficiently trained 
without extra costs for boundary and interface treatments. 
In addition, we note that this kind of penalty-free boundary and interface treatments  is 
proposed in  \cite{sheng2022pfnn2} for fully connected neural networks with residual connections, but the novelty of our F-D3M
lies on the Fourier feature networks for high frequency problems.

The rest of the paper is organized as follows. 
Section~\ref{sec:pdes} gives standard procedures of deep learning methods for solving PDEs.  
In section~\ref{sec:algo}, we first discuss standard overlapping domain decomposition methods
and boundary treatments for the neural network based methods, and next present details of F-D3M 
algorithm.  In section~\ref{sec:experiments}, 
we demonstrate the efficiency of F-D3M with numerical experiments, where both diffusion and Helmholtz 
problems are studied. Finally  section~\ref{sec:conclusions} concludes the paper. 

\section{Deep learning for PDEs} \label{sec:pdes}
Let $\Omega \subset \mR^d$ ($d \in \mN$) denote a spatial domain which is  bounded, connected and with a
polygonal boundary $\pOmega$, 
and $\bx\in \mR^{d}$ denote a spatial variable.
The physics of problems considered 
are governed by 
a PDE over the domain $\Omega$ and
boundary conditions on the boundary $\pOmega$, which are stated as: 
find $u(\bx)$ mapping $\Omega$ to $\mR$, such that
\begin{eqnarray}
    \pL\left(u\left(\bx\right)\right)=f(\bx) \quad &\forall\bx\in \Omega, \label{pde1}\\
    \pB\left(u\left(\bx\right)\right)=g(\bx) \quad &\forall\bx\in \pOmega,\label{pde2}
\end{eqnarray}
where $\pL$ is a partial differential operator, $\pB$ is a boundary operator, $f(\bx)$ is the source function and $g(\bx)$ specifies the boundary conditions. 

When applying deep learning methods for \eqref{pde1}--\eqref{pde2}, 
the target solution $u(\bx)$ is approximated by a neural network $u(\bx;\btheta)$ where $\btheta$ represents the parameters of the neural network. 
Defining $||u(\bx)||_{2,\Omega}^{2}:=\int_{\Omega}|u(\bx)|^2d\bx$ and $||u(\bx)||_{2,\pOmega}^{2}:=\int_{\pOmega}|u(\bx)|^2d\bx$, 
the neural network approximation $u(\bx;\btheta)$ is obtained by finding the optimal $\btheta$, 
which minimizes the following composite loss function
\begin{equation}
        J\left(u\left(\bx;\btheta\right)\right)=J_r(u(\bx;\btheta))+\lambda J_b(u(\bx;\btheta)), \\
    \label{com_loss_func}
\end{equation}
where
\begin{equation}
    \begin{aligned}
        J_r\left(u\left(\bx;\btheta\right)\right)&:=\big\|\pL\left(u\left(\bx;\btheta\right)\right)-f(\bx)\big\|_{2,\Omega}^2=\int_{\Omega}\big|\pL\left(u\left(\bx;\btheta\right)\right)-f(\bx)\big|^2d\bx,\\
        J_b\left(u\left(\bx;\btheta\right)\right)&:=\big\|\pB\left(u\left(\bx;\btheta\right)\right)-g(\bx)\big\|_{2,\pOmega}^2=\int_{\pOmega}\big|\pB\left(u\left(\bx;\btheta\right)\right)-g(\bx)\big|^2d\bx,
    \end{aligned}
    \label{Jr_and_Jb}
\end{equation}
and $\lambda>0$ is a weight coefficient to be determined. 
In \eqref{Jr_and_Jb}, 
$J_r(u(\bx;\btheta))$ measures how well $u(\bx;\btheta)$ satisfies the governing  equation \eqref{pde1} and $J_b(u(\bx;\btheta))$ is used to contain $u(\bx;\btheta)$ to meet the boundary condition \eqref{pde2}. 

In practice, the integrals in \eqref{Jr_and_Jb}  are typically approximated through numerical integration.
That is,  choosing two sets of randomly sampled collocation points $\mS_r:=\{\bx_r^{(n)}\}_{n=1}^{N_r}\subset \Omega$ and $\mS_b:=\{\bx_b^{(n)}\}_{n=1}^{N_b}\subset \pOmega$, 
and denoting $\|v(\bx)\|_{\mS_r}^2:=\frac{1}{N_r}\sum_{n=1}^{N_r}|v(\bx_r^{(n)})|^2$ and $\|v(\bx)\|_{\mS_b}^2=\frac{1}{N_b}\sum_{n=1}^{N_b}|v(\bx_b^{(n)})|^2$ for any arbitrary function $v$,  
the following empirical loss function is defined
\begin{equation}
    \begin{aligned}
        J_{E}\left(u\left(\bx;\btheta\right)\right)=&\big\|\pL\left(u\left(\bx;\btheta\right)\right)-f(\bx)\big\|_{\mS_r}^2+\lambda \big\|\pB\left(u\left(\bx;\btheta\right)\right)-g(\bx)\big\|_{\mS_b}^2\\
        =&\frac{1}{N_r}\sum_{n=1}^{N_r}\big|\pL\left(u\left(\bx_r^{(n)};\btheta\right)\right)-f\left(\bx_r^{(n)}\right)\big|^2\\
        &+\frac{\lambda}{N_b}\sum_{n=1}^{N_b}\big|\pB\left(u\left(\bx_b^{(n)};\btheta\right)\right)-g\left(\bx_b^{(n)}\right)\big|^2,\\
    \end{aligned}
    \label{em_loss_func}
\end{equation}
based on which we choose the optimal parameter $\btheta^*$:
\begin{equation}
        \btheta^*=\arg \min_{\btheta}J_{E}\left(u\left(\bx;\btheta\right)\right).
    \label{optimal_theta}
\end{equation}
It should be noted that values of $\pL(u(\bx_r^{(n)};\btheta))$ (for $n=1,\ldots,N_r$) in \eqref{em_loss_func} 
can be easily computed by the automatic differentiation for neural networks \cite{baydin2018automatic}.  
For the optimization problem \eqref{optimal_theta}, stochastic gradient descent methods \cite{bottou2018optimization} are usually used. As discussed in section \ref{sec:intro}, this kind of deep learning based methods for PDEs 
has the problem caused by the  frequency principle, which leads to difficulties for solving problems with high frequency
modes. For this purpose, this paper develops a new domain decomposition based deep learning approach for these 
challenging situations, of which details are presented in the next section. 

\section{A Fourier feature based deep domain decomposition method} \label{sec:algo} 
In this section, we first review the related domain decomposition methods, and then discuss efficient boundary 
treatments for neural network approximations. After that, a new multi Fourier feature network (MFFNet) is 
presented, which is used to construct approximations for local solutions posed in subdomains. 
Finally, details of our Fourier feature based deep domain decomposition method (F-D3M) 
and the corresponding algorithm  are presented.

\subsection{Overlapping domain decomposition methods}
We review overlapping domain decomposition methods following the presentation in \cite{quarteroni1999domain,toselli2004domain}. 
The idea of overlapping domain decomposition is originally  introduced by Schwarz \cite{schwarz1870ueber}  
 and its parallel version is further developed by Lions \cite{lions1989schwarz}, 
 which is referred to as additive Schwarz methods. 
 Geometry convergence results of the Schwarz alternating procedure for elliptic PDEs are given 
 by Chan et al.\ \cite{chan1991geometry}. 
 Specific overlapping domain decomposition methods are 
 actively developed and  widely used for many types of PDEs,
 and some related works are listed as follows. 
 Domain decomposition algorithms for indefinite elliptic problems are developed by Cai and Widlund \cite{cai1992domain},
 and methods for  parabolic problems  are developed in \cite{cai1991additive,cai1994multiplicative,li2018multilevel,deng2022efficient}.
Cai et al.\ present overlapping domain decomposition methods for the Helmholtz equation \cite{cai1998overlapping},
and Chen et al.\ present a robust method for the Helmholtz equation with high wave numbers \cite{chen2016robust}. 
An overlapping domain decomposition preconditioner is proposed  for discontinuous Galerkin approximations of advection-diffusion problems by Lasser and Toselli \cite{lasser2003overlapping},
and a balancing domain decomposition method by constraints of advection-diffusion problems
is proposed by Tu and Li \cite{tu2008balancing}.  
Multi-scale domain decomposition methods are established by Aarnes and Hou \cite{aarnes2002multiscale}
and overlapping domain decomposition preconditioners  for multiscale flows  are established
by Galvis and Efendiev \cite{galvis2010domain}. 

In overlapping domain decomposition methods, the domain $\Omega$ is decomposed into  
$N_{d}$ ($N_d> 1$) overlapping subdomains such that 
\begin{equation}
        \Omega=\Omega_1 \cup \Omega_2 \cup \cdots \cup \Omega_{Nd}.
    \label{domain_decomposition}
\end{equation}
\begin{figure}[!ht]
    \centerline{
    \begin{tabular}{cc}
    \includegraphics[width=8cm]{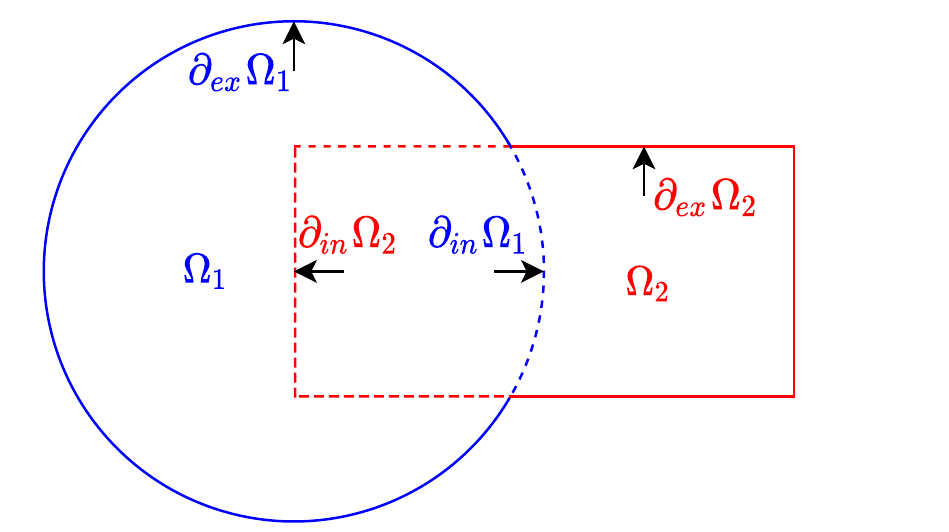}
    \end{tabular}}
    \caption{Example of a two-subdomain system with overlapping domain decomposition: $\Omega=\Omega_1 \cup \Omega_2$, 
    $\partial_{in}\Omega_1=\partial_{2}\Omega_1=\pOmega_1\cap\Omega_2$, 
    $\partial_{ex}\Omega_1=\pOmega_1\cap\pOmega$,
    $\partial_{in}\Omega_2=\partial_{1}\Omega_2=\pOmega_2\cap\Omega_1$ and
    $\partial_{ex}\Omega_2=\pOmega_2\cap\pOmega$.
    }
    \label{fig_dd_example}
\end{figure}
For each subdomain $\Omega_i$ ($i=1,\ldots,N_d$), $\pOmega_i$ denotes the set of its boundaries, and $\Lambda_i$ represents the set of its neighboring subdomain indices, i.e., $\Lambda_i:=\{j\,|\,j\in\{1,\dots,N_{d}\},j \neq i$ and $\Omega_i \cap \Omega_j \neq \emptyset \}$. Each boundary set $\pOmega_i$ is splited into two  parts: $\pOmega_i=\pexoi \cup \partial_{in}\Omega_i$, where $\pexoi:=\pOmega_i\cap \pOmega$ are exterior boundaries and $\partial_{ in} \Omega_{i}:=\cup_{j\in \linteri}\{\partial_{j} \Omega_i\}$ where $\partial_j\Omega_{i}:=\pOmega_i\cap\Omega_j$ are interfaces introduced by domain decomposition. 
A notional example demonstrating this notation for two subdomains with single interface is shown in Figure~\ref{fig_dd_example}. To formulate local problems, we introduce the decomposited local operators $\{\pL_i:=\pL\,|_{\Omega_i}\}_{i=1}^{N_d}$ and $\{\pB_i:=\pB\,|_{\pexoi}\}_{i=1}^{N_d}$ and local functions $\{f_i(\bx):=f(\bx)\,|_{\Omega_i}\}_{i=1}^{N_d}$ and $\{g_i(\bx):=g(\bx)\,|_{\pexoi}\}_{i=1}^{N_d}$. Let $\{\phi_{i,j}(\bx)\}_{j\in\Lambda_i}$ denote the interface functions defined on $\{\pjoi\}_{j\in\Lambda_i}$.  
The overlapping domain decomposition methods proceed as follows. Given initial guess $\phi_{i,j}^0(\bx)$
for each  interface function,  
for each iteration step $k=0,1,\dots$, we solve $N_d$ local problems: find $u_i(\bx)$ mapping
$\Omega_i$ to $\mR$ for $i=1,2,\dots,N_d$, such that
\begin{eqnarray}
        \pL_i\left(u_i^{k}\left(\bx\right)\right) & = & f_i(\bx) \quad \forall\bx\in \Omega_i,\label{dd1} \\ 
        \pB_i\left(u_i^{k}\left(\bx\right)\right) & = & g_i(\bx) \quad \forall\bx\in \pexoi,\label{dd2}\\
       u_i^{k}\left(\bx\right) & = & \phi_{i,j}^{k}(\bx) \quad \forall\bx\in \pjoi,\,\text{for all}\,j\in\Lambda_i,\label{dd3}
\end{eqnarray}
where $\phi_{i,j}^{k}(\bx)$ is the interface function at the $k$-th iteration step and is updated based on the local solutions on the neighbouring subdomains through
\begin{equation}
    \phi_{i,j}^{k+1}(\bx):=u_j^{k}(\bx) \quad \forall \bx \in \pjoi,\, \text{for all}\,j\in\Lambda_i.
    \label{boundary_function_update}
\end{equation}

This paper focuses on applying overlapping domain decomposition to construct neural network 
approximations with \eqref{em_loss_func}. Our motivation is based on the fact that frequencies involved in 
the underlaying PDE problems are typically defined relative to the size of the spatial domain. 
Even though there may exist high frequencies in the original global problem, 
through dividing the global spatial domain into small subdomains, the relative frequencies posed on
subdomains can become smaller. Then, neural network approximations for local solutions can be efficiently 
trained, and details of our full approach are presented at the end of this section.  

\subsection{Interface treatments}
In our original D3M setting \cite{li2019d3m}, the boundary condition for interfaces is treated through separating
each local solution into a particular function to satisfy \eqref{dd2}--\eqref{dd3} and an auxiliary solution to satisfy homogeneous
boundary conditions, and deep neural networks are constructed to approximate the  auxiliary solution. 
To meet the homogeneous boundary condition (i.e., let $\phi_{i,j}^{k}(\bx)=0$ in \eqref{dd3}),
a loss function with the norm of the auxiliary solution restricted on interfaces is introduced in \cite{li2019d3m}. 
As finding the optimal weights to balance the terms for boundary  conditions and the governing PDEs
remains an open challenging problem,    
we  apply the boundary treatment method proposed in \cite{berg2018unified,sheng2021pfnn,sheng2022pfnn2} as follows.

For each local problem \eqref{dd1}--\eqref{dd3}, a smooth particular function $\Phi_i^k(\bx)$ (for $i=1,\ldots,N_d$)
is introduced to satisfy the boundary conditions, i.e., 
\begin{eqnarray}
    \pB_i\left(\Phi_i^{k}(\bx)\right)&=&g_i(\bx)\quad \forall \bx \in \pexoi, \label{p1}\\
    \Phi_i^{k}(\bx)&=&\phi_{i,j}^{k}(\bx)\quad \forall\bx\in \pjoi,\, \text{for all}\,j\in\Lambda_i. \label{p2}
\end{eqnarray}
Then,  a neural network model to approximate each local solution $u_i^{k}(\bx) $ of \eqref{dd1}--\eqref{dd3} is 
written as
\begin{equation}
    \begin{aligned}
    u_i^{k}(\bx;\btheta_i):=\Phi_i^{k}(\bx)+D_i(\bx)\hat u_i^{k}(\bx;\btheta_i),
    \end{aligned}
    \label{anstz}
\end{equation}
where $\hat u_i^k(\bx;\btheta_i)$ represents a neural network with parameters $\btheta_i$  
and $D_i(\bx)$ is a smooth distance function defined in $\Omega_i$, which satisfies 
\begin{eqnarray}
    D_i\left(\bx\right)&=&0 \quad \forall \bx \in \pOmega_i,  \label{Dx}
\end{eqnarray}
and increases away from the boundary. 
For the problems with simple spatial domains, 
$\Phi_i^{k}(\bx)$ and $D_i(\bx)$ can typically  be derived analytically. 
When the spatial domain is complex, 
neural network approximations for  $\Phi_i^{k}(\bx)$ and $D_i(\bx)$ are recommended \cite{berg2018unified,sheng2021pfnn}. 

In our setting, $u_i^k(\bx;\btheta_i)$ satisfies the conditions \eqref{dd2}--\eqref{dd3} automatically. 
To simplify the presentation, for any arbitrary set $\mS$, $|\mS|$ denotes its size. 
With a set of randomly sampled collocation points $\mS_i:=\{\bx_{n}\}_{n=1}^{|\mS_i|}\subset \Omega_i$, 
the following loss function is defined for each local problem, 
\begin{equation}
    J_{E_i}^{k}\left(u^{k}_i\left(\bx;\btheta_i\right)\right)=\frac{1}{|\mS_{i}|}\sum_{\bx_{n}\in\mS_i}\left|\pL_i\left(u^{k}_i\left(\bx_{n};\btheta_i\right)\right)-f_i\left(\bx_{n}\right)\right|^2 ,\\
    \label{uncon_opti_problem}
\end{equation}
and we choose the optimal parameter $\btheta_i^*$:
\begin{equation}
    \btheta_i^{*}=\arg \min_{\btheta_i}J_{E_i}^{k}\left(u^{k}_i\left(\bx;\btheta_i\right)\right).
    \label{opti_btheta_i}
\end{equation}
As discussed in section \ref{sec:pdes}, the stochastic gradient descent (SGD) method is typically used to solve
this kind of optimization problems \cite{bottou2018optimization}.  
For the loss function \eqref{uncon_opti_problem}, 
with an initial guess  $\btheta^{(0)}_i$,
the parameters $\btheta^{(s)}_i$ at $s$-th iteration step of SGD are updated as,  
\begin{equation}
    \begin{aligned}
        \btheta_i^{(s+1)}=\btheta_i^{(s)}-\alpha\nabla_{\btheta}\left[\frac{1}{|\mS_{i}|}\sum_{\bx_{n}\in\mS_i}\left|\pL_i\left(u^{k}_i\left(\bx_{n};\btheta_i^{(s)}\right)\right)-f_i\left(\bx_{n}\right)\right|^2\right],
    \end{aligned}
    \label{SGD}
\end{equation}
where $\alpha$ is a learning rate and  $s=0,1,\ldots$. 
The set of collocation points $\mS_i$ is generated by some given distribution at each  iteration step. 
In this work, the Adam optimizer version of SGD is employed, 
which adopts adaptive learning rates ($\alpha$ in \eqref{SGD}) for different components of the parameters through estimating the first and the second moments of the gradients \cite{kingma2014adam}.

\subsection{Multi Fourier Feature Network (MFFNet)}
\label{section_MFFNet}
As discussed in \cite{xu2020frequency,luo2019theory}, due to the  frequency principle, there exists difficulties for  neural networks to approximate
high frequency modes. 
To relieve this limitation, random Fourier feature (RFF) mappings are developed  in \cite{tancik2020fourier}
and multi random Fourier feature mappings  are developed in \cite{wang2021eigenvector}. 
In addition, deterministic Fourier feature mapping (also called position encoding) are proposed to fit data with high frequency variation \cite{mildenhall2020nerf},  where the inputs of the networks include the original independent variables and the mapped variables.  
Inspired by \cite{mildenhall2020nerf}, to balance the ability for learning low frequencies 
and high frequencies, a variant random Fourier feature mapping that keeps original independent variables in the input, is proposed in this work, i.e.,
\begin{equation}
    \bg(\bx)=\big[\cos(2\pi\bB\bx),\bx,\sin(2\pi\bB\bx)\big]^{\text{T}},
    \label{vRFF}
\end{equation}
where  each entry of $\bB\in \mR^{m\times d}$ ($m$ is a given integer) is sampled from an isotropic Gaussian distribution $\mathcal{N}(0,\sigma^2)$ and the standard deviation  $\sigma$ is a hyper-parameter. 
\begin{figure}[!ht]
    \centerline{
    \begin{tabular}{cc}
    \includegraphics[width=14cm]{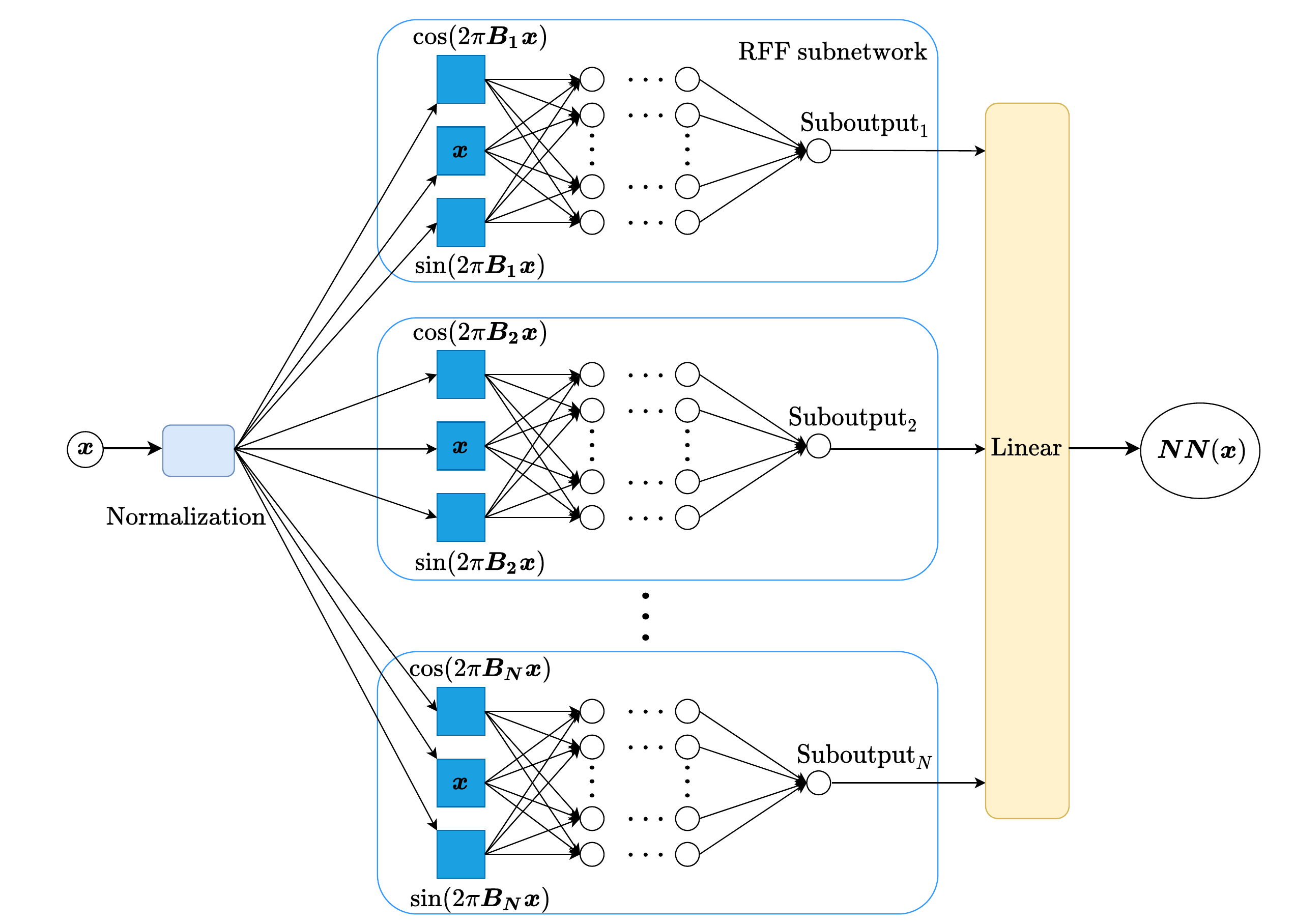}
    \end{tabular}}
    \caption{Multi Fourier Feature Network (MFFNet).}
    \label{fig_MRFFnet}
\end{figure}

To efficiently capture modes with different  frequencies,  
we design a new multi Fourier feature network (MFFNet), which is shown in Figure~\ref{fig_MRFFnet}.
A MFFNet consists of $N$ subnetworks with different parameters $\bB_n$ ($n=1,\ldots,N$),
and each subnetwork is called a RFF subnetwork. The overall output of MFFNet is a weighted sum of the outputs of subnetworks. In each RFF subnetwork, the input variable is first normalized by a linear transformation, and
 the corresponding random Fourier feature is computed through \eqref{vRFF}. 
 The random Fourier feature
 for each RFF subnetwork is denoted by $\bg_n(\bx)=[\cos(2\pi \bB_n \bx),\bx,\sin(2\pi\bB_n\bx)]^{\text{T}}$, 
 where the parameters $\bB_n$ determine the preferred frequencies of the subnetwork.  
 Then $\bg_n(\bx)$ goes through a fully connected architecture to produce a subnetwork output $\boldsymbol{H}_n(\bx)$. Finally, All the subnetwork outputs are concatenated with a linear layer to give the final output $\boldsymbol{NN}(\bx)$.  
The detailed workflow of the MFFNet is defined as
\begin{eqnarray}
        \bg_n(\bx)&=&\big[\cos(2\pi\bB_n\bx),\bx,\sin(2\pi\bB_n\bx)\big]^{\text{T}}, \quad n=1,2,...,N, \\
        \boldsymbol{H}_n(\bx)&=&\mathcal{FCN}_n(\bg_n(\bx)), \quad n=1,2,...,N, \\
        \boldsymbol{NN}(\bx)&=&\boldsymbol{W}\cdot\big[\boldsymbol{H}_1,\boldsymbol{H}_2,\cdots,\boldsymbol{H}_N\big]+\boldsymbol{b},
\end{eqnarray}
where $\mathcal{FCN}_n$ represents a fully connected network architecture, and $\boldsymbol{W}$ and $\boldsymbol{b}$ denote the weights and the bias of the last linear layer.  
It should be noted that different subnetworks in MFFNet are independent 
and their sizes  can be separately adjusted.

\subsection{F-D3M algorithm} 
Given a PDE problem \eqref{pde1}--\eqref{pde2},  
we decompose the spatial domain $\Omega$ into $N_d$ overlapping subdomains such that 
$\Omega=\cup_{i=1}^{N_d}\Omega_i$, and our F-D3M algorithm proceeds as follows. 
To start with, for each subdomain $\Omega_i$, 
a distance function $D_i(\bx)$ according to $\pOmega_i$ is constructed (see \eqref{Dx}), 
and  a MFFNet is initialized (see section \ref{section_MFFNet}),  which is denoted by  $\hat u_i^0(\bx;\btheta_i)$. 
Then, we set initial iteration step to $k=0$ and denote $u_i^0(\bx;\btheta_i):=\hat u_i^0(\bx;\btheta_i)$. 

As each local solution consists of two parts---the boundary part and the interior part shown in \eqref{anstz}, 
the main procedure of F-D3M consists of two alternating iterations. 
The first one is to update the particular function $\Phi_i^{k+1}(\bx)$ for boundary conditions; the second one is to solve the local problem \eqref{uncon_opti_problem}--\eqref{opti_btheta_i}. 
The particular function $\Phi_i^{k+1}(\bx)$ stratifies the exterior boundary condition \eqref{p1} 
and the interface condition \eqref{p2}. As the exterior boundary condition is fixed, only 
the interface condition  at each iteration step needs to be updated,
and based on \eqref{boundary_function_update}, this updating procedure is conducted through
\begin{eqnarray}
\phi^{k+1}_{i,j}(\bx):=\left.u_j^{k}(\bx;\btheta_j)\right|_{\bx\in\pjoi}, 
\end{eqnarray}
where $j\in\Lambda_i$ and $\Lambda_i$ includes the indices of the neighboring  subdomains of  $\Omega_i$.  
After the particular function  $\Phi_i^{k+1}(\bx)$  is updated, we train a deep neural network for the interior part of the local solution (i.e., $u_i^{k+1}(\bx;\btheta_i)$ in  \eqref{anstz}).

Once all local problems are solved, a global solution can be assembled.  
Letting $\alpha_{\bx}$ denote the set consisting of indices of subdomains that contain $\bx$, i.e., $\alpha_{\bx}:=\{i\,|\,i\in\{1,\dots,N_{d}\},\bx \in \Omega_i\}$ and $|\alpha_{\bx}|$ is the size of $\alpha_{\bx}$, 
a global solution $u^{k+1}(\bx;\btheta)$ at iteration step $k+1$ is assembled as 
\begin{equation}
    u^{k+1}(\bx;\btheta)=\frac{1}{|\alpha_{\bx}|}\sum_{s\in\alpha_{\bx}}u_s^{k+1}\left(\bx;\btheta_s\right).
    \label{assemble}
\end{equation}
We assess the difference between $u^{k+1}(\bx;\btheta)$ and $u^{k}(\bx;\btheta)$ 
using the quantity $\epsilon(u^{k+1}(\bx;\btheta),u^{k}(\bx;\btheta))$---$\epsilon(\cdot,\cdot)$ is defined as
\begin{equation}
    \begin{aligned}
        \epsilon\left(u_1(\bx),u_2(\bx)\right)=\frac{\left(\sum_{n=1}^{N_{t}}\left|u_1\left(\bx_t^{(n)}\right)-u_2\left(\bx_t^{(n)}\right)\right|^2\right)^{\frac{1}{2}}}{\left(\sum_{n=1}^{N_{t}}\left|u_2\left(\bx_t^{(n)}\right)\right|^2\right)^{\frac{1}{2}}},
    \end{aligned}
    \label{relative_l2_error}
\end{equation}
where $u_1(\bx)$ and $u_2(\bx)$ are two arbitrary functions, and $\{\bx_t^{(n)}\}_{n=1}^{N_{t}}$ are $N_{t}$ uniformly sampled test points with range $\Omega$. 
After that, we update the iteration step through setting $k:=k+1$, 
and repeat the above procedures until the difference  $\epsilon(u^{k+1}(\bx;\btheta),u^{k}(\bx;\btheta))$ 
is smaller than a given tolerance $tol$. 

Details of the F-D3M algorithm are summerized in Algorithm~\ref{alg_F_D3M}. 
Here, $\eta$ generically denotes the difference between the approximation solutions at two adjacent 
iteration steps (line  24 of Algorithm~\ref{alg_F_D3M}) and it is initialized with an arbitrary value such that $\eta>tol$.
The number of epochs for SGD (see \eqref{SGD}) is denoted by $N_e$.  
At each iteration step of SGD, collocation points are generated with a given distribution, 
and the uniform distribution is applied in this work.
Cf.\ \cite{tang2022adaptive} for systematic generative models for generating collocation points. 
\begin{algorithm}
    \caption{F-D3M for PDEs}
    \label{alg_F_D3M}
    \begin{algorithmic}[1]
    \State \textbf{Input:}  The PDE problem, the global $\Omega$ and the overlapping subdomains $\left\{\Omega_i\right\}_{i=1}^{N_d}$.     
    \State Initialize MFFNets $\hat u_i^{0}\left(\bx;\btheta_i\right)$ (see Figure~\ref{fig_MRFFnet}),
               and set  $u_i^{0}\left(\bx;\btheta_i\right)=\hat u_i^{0}\left(\bx;\btheta_i\right)$ for   $i=1,2,\cdots,N_d$.
    \State Set $k=0$, and initialize an arbitrary value for $\eta$ such that $\eta> tol$.  
    \For{$i=1:N_d$}
    \State Construct $D_i\left(\bx\right)$ using the boundary coordinates of $\Omega_i$ (see \eqref{Dx}).
    \EndFor
    \While{$\eta>tol$}
    \For{$i=1:N_d$}
        \For{$j\in\Lambda_i$}
            \State Let $ \phi^{k+1}_{i,j}(\bx):=\left.u_j^{k}(\bx;\btheta_j)\right|_{\bx\in\pjoi}$. 
        \EndFor
        \State Construct $\Phi_i^{k+1}\left(\bx\right)$ using $g_{i}\left(\bx\right)$ and $\phi^{k+1}_{i,j}(\bx)$ (see \eqref{p1}--\eqref{p2}).
    \EndFor
    \For{$i=1:N_d$}
    \State Let $\btheta_i^{(0)}:=\btheta_{i}$.
    \State Let $u_i^{k+1}\left(\bx;\btheta_i^{(0)}\right)=\Phi_i^{k+1}\left(\bx\right)+D_i\left(\bx\right)\hat u_i^{k+1}\left(\bx;\btheta_i^{(0)}\right)$.
    \For {$s=0:N_e-1$}
        \State Generate samples $\mS_{i}=\left\{\bx_{n}\right\}_{n=1}^{|\mS_i|}$ by a given distribution with range $\Omega_i$.
        \State Update $\btheta_i^{(s)}$ using the stochastic gradient descent method (see \eqref{SGD}).
     \EndFor
    \State Let $\btheta_i:=\btheta_i^{(N_e-1)}$.
    \EndFor
    \State Assemble the global solution $u^{k+1}\left(\bx;\btheta\right)$ with $\left\{u_i^{k+1}\left(\bx;\btheta_i\right)\right\}_{i=1}^{N_d}$ (see \eqref{assemble}).
    \State Compute $\eta=\epsilon\left(u^{k+1}\left(\bx;\btheta\right),u^{k}\left(\bx;\btheta\right)\right)$ (see \eqref{relative_l2_error}).
    \State Let $k=k+1$.
    \EndWhile
    \State Set $u\left(\bx;\btheta\right)=u^{k}\left(\bx;\btheta\right)$.
    \State \textbf{Output:} The overall F-D3M solution $u\left(\bx;\btheta\right)$. 
    \end{algorithmic}
\end{algorithm}

\section{Numerical experiments} \label{sec:experiments}
In this section, 
numerical experiments are conducted to illustrate the effectiveness of our F-D3M approach presented in Algorithm
\ref{alg_F_D3M}. Three test problems are considered---the first one is a one-dimensional diffusion problem, 
the second one is a two-dimensional diffusion problem, and the third one is a Helmholtz equation with a high 
frequency. 
\subsection{One-dimensional diffusion problem} 
We start with the following one-dimensional Poisson equation, 
\begin{equation}
    \begin{aligned}
        \Delta u(x)&=f(x)\quad \text{in}\ \Omega, \\
        u(x)&=g(x)\quad \text{on}\ \pOmega, \\
    \end{aligned}
    \label{test1-pde}
\end{equation}
where the spatial domain is $\Omega=[-1,1]$ and the exact solution is 
\begin{equation}
    \begin{aligned}
        u(x)=\sin(5\pi x)+\sin(30\pi x).     \label{test1-sol}
    \end{aligned}
\end{equation}
In \eqref{test1-pde}, $f$ and $g$ are specified with \eqref{test1-sol}.
It is clear that the solution \eqref{test1-sol} contains two high frequency components. 

For this test problem, the spatial domain $\Omega$ is decomposed into $N_d=5$ overlapping subdomains
$\{\Omega_i\}_{i=1}^{N_d}$  and   
each subdomain is defined as $\Omega_i=[x_{left}^{(i)},x_{right}^{(i)}]$ with    
\begin{equation}
    \begin{cases}
        x_{left}^{(i)}&=\max\left\{-1,-1+(i-1)\times\frac{2}{N_d}-\frac{w}{2}\right\},\\
        x_{right}^{(i)}&=\min\left\{1,-1+i\times\frac{2}{N_d}+\frac{w}{2}\right\},
    \end{cases}
    \label{x1}
\end{equation}
where the width of each overlapping region is set to $w=0.2$. 
In our F-D3M algorithm (Algorithm \ref{alg_F_D3M}), a MFFNet (see Figure \eqref{fig_MRFFnet}) with two RFF subnetworks is used for each subdomain. 
The standard deviation $\sigma$ is set to $1$ and $30$ for the first and the second subnetworks respectively, 
and each subnetwork has one hidden layer with ten neurons.  
For comparison, two methods for solving this test problem in the global domain without domain decomposition 
are considered---one is using the MFFNet for the whole domain, 
and  the other is using a conventional fully connected network (FCN). 
These two methods are referred to as Global-MFFNet and Global-FCN in the following.
For  both methods,  the weight coefficient in \eqref{em_loss_func} is set to $\lambda=100$. 
The setting of MFFNet in Global-MFFNet is the same as that for F-D3M, while Global-FCN has
one hidden layer with twenty neurons. 
For all networks, the hyperbolic tangent function is employed as the activation function 
and the Kaiming scheme is used for  initialization \cite{he2015delving} . 
All networks are trained by the Adam \cite{kingma2014adam} optimizer with defaulting settings. 
An exponential learning rate with the initial learning rate of 0.01 and a decay rate of 0.9 every 1000 training epochs is employed. 
The number of epochs for each domain decomposition iteration in F-D3M is set to $2500$ and the number of iterations is set to $20$. For each training epoch, $|\mS_i|=400$ collocation points are generated 
using the uniform distribution with range $\Omega_i$ for $i=1,\ldots,N_d$. 
To train Global-FCN and Global-MFFNet, $50000$ epochs are used and 
the number of collocation points for the interior domain and the boundary 
are set to $N_r=2000$ and $N_b=2$ (the boundary consists of  the two endpoints for this one dimensional problem) respectively.
\begin{figure}[!ht]
    \centerline{
    \begin{tabular}{ccc}
    \includegraphics[width=6cm]{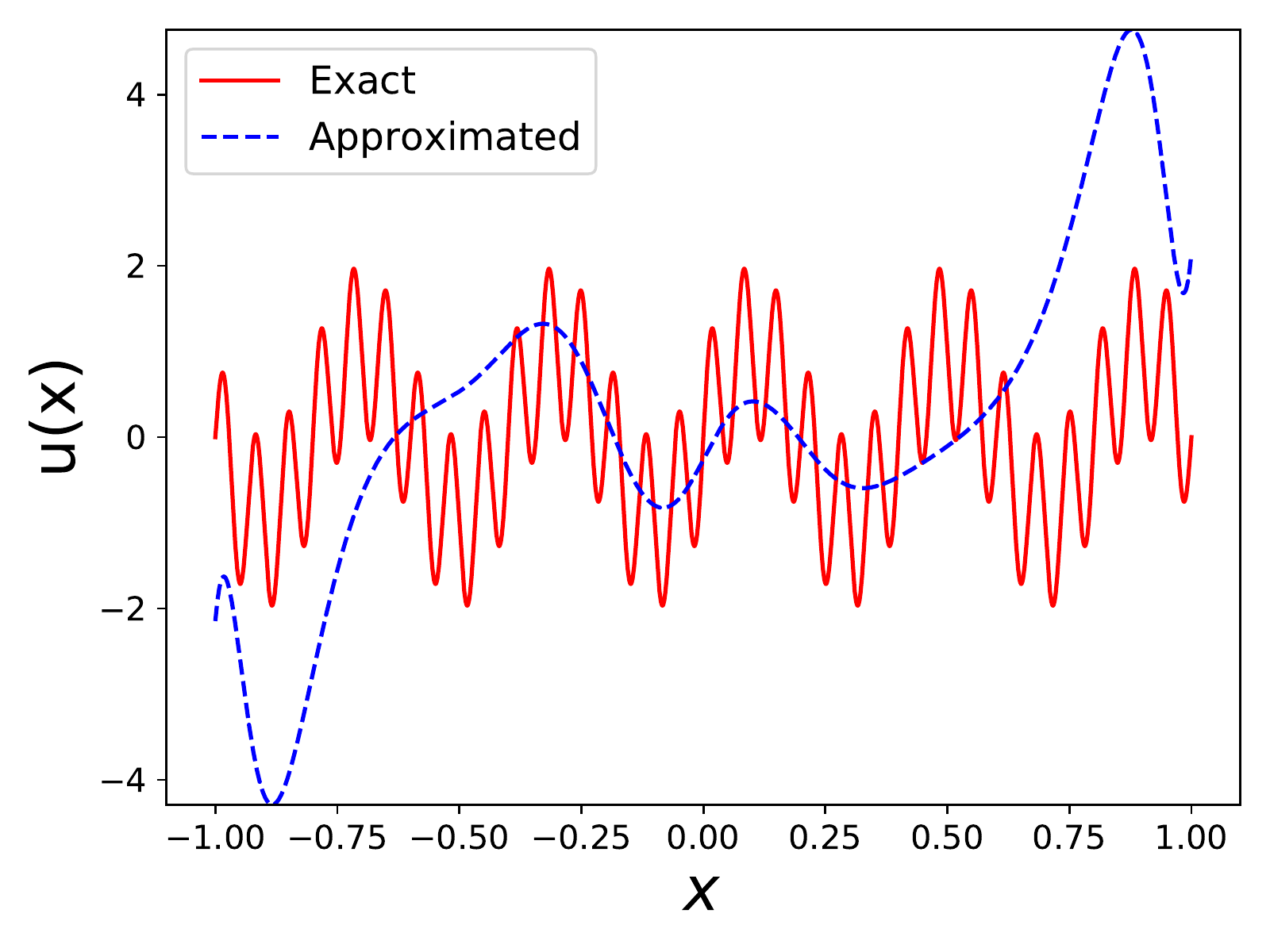}&
    \includegraphics[width=6cm]{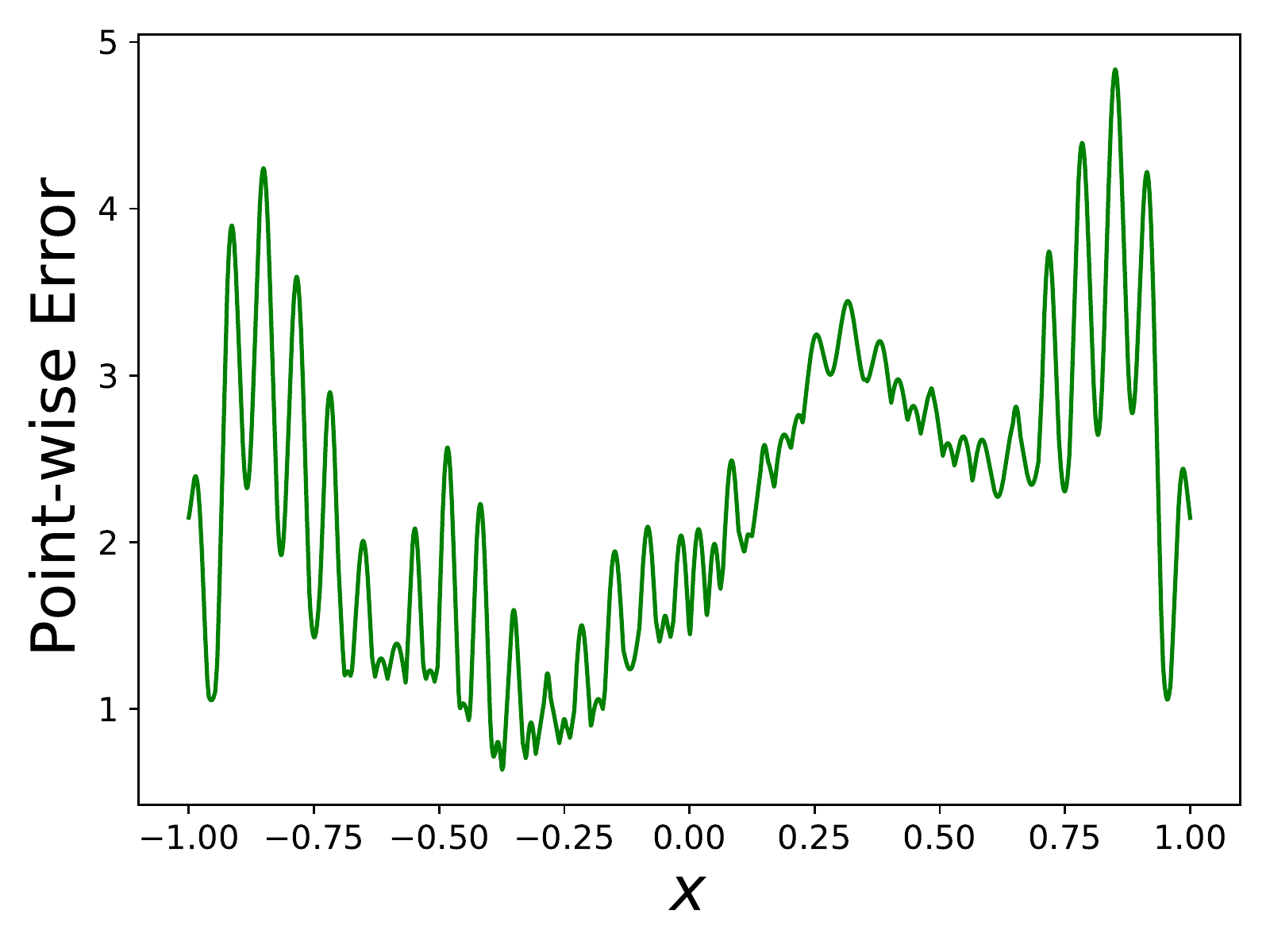}\\
    (a) Global-FCN & (b) Global-FCN point-wise error \\
    \includegraphics[width=6cm]{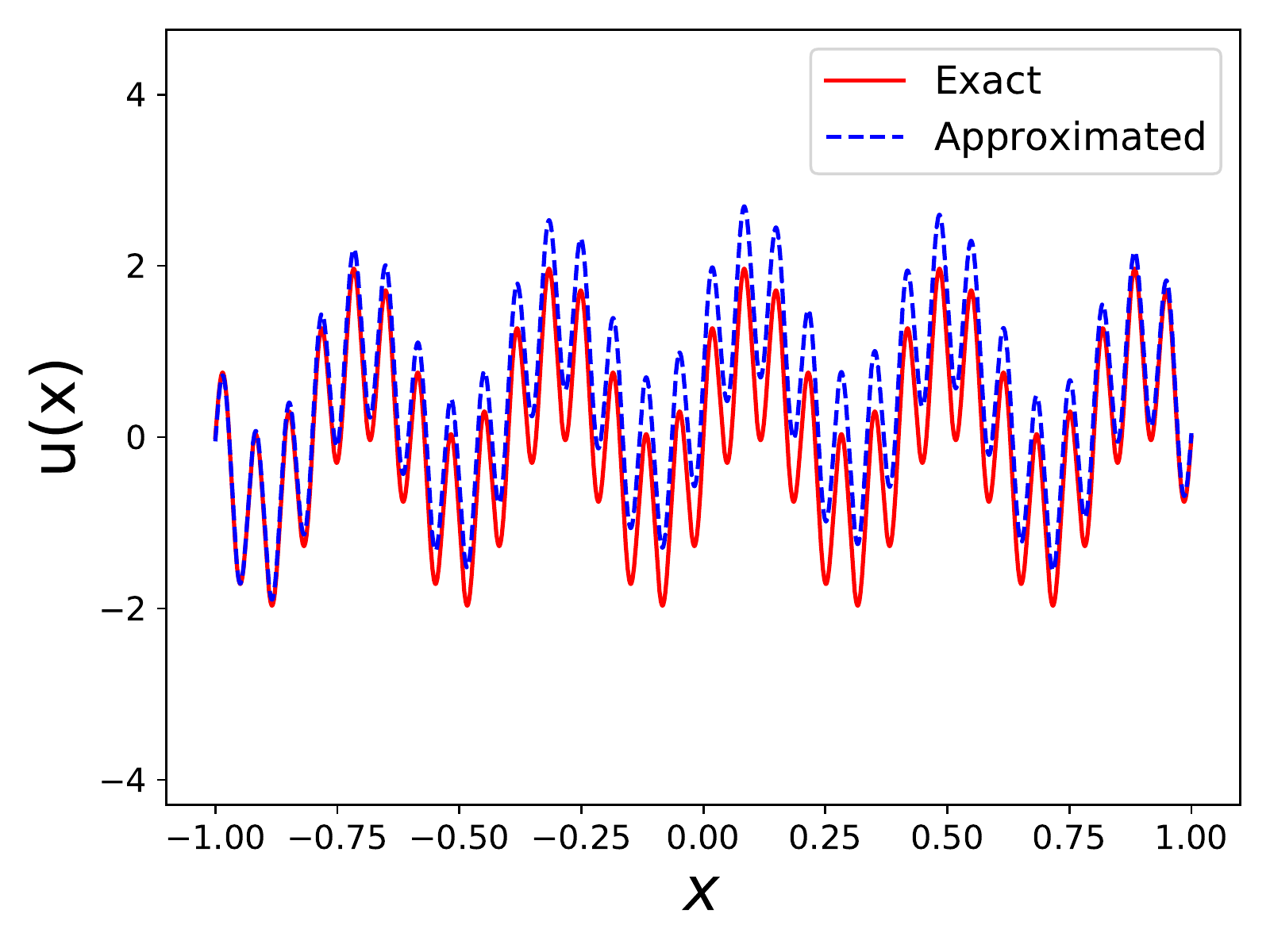}&
    \includegraphics[width=6cm]{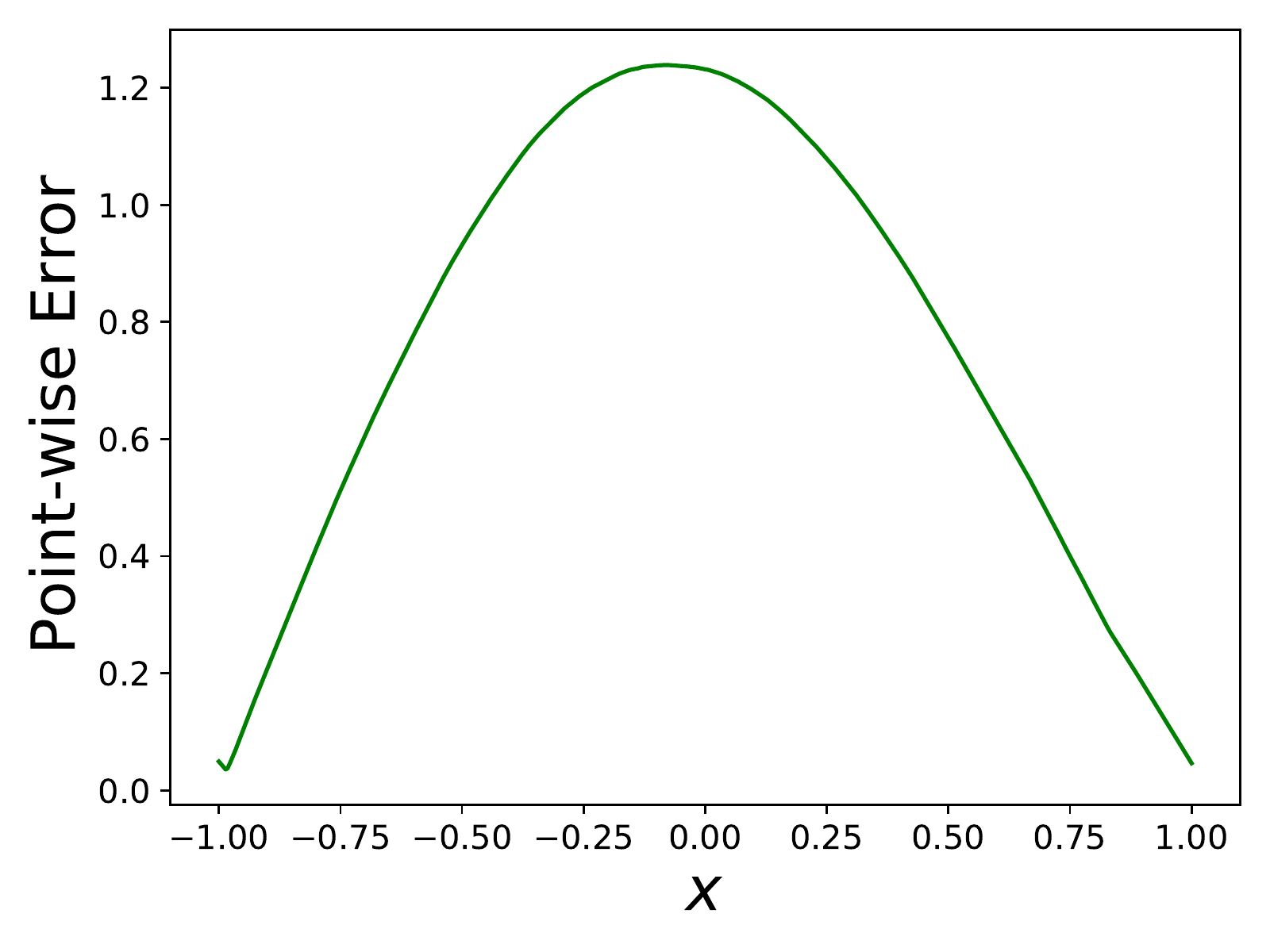}\\
    (c) Global-MFFNet & (d) Global-MFFNet point-wise error \\
    \includegraphics[width=6cm]{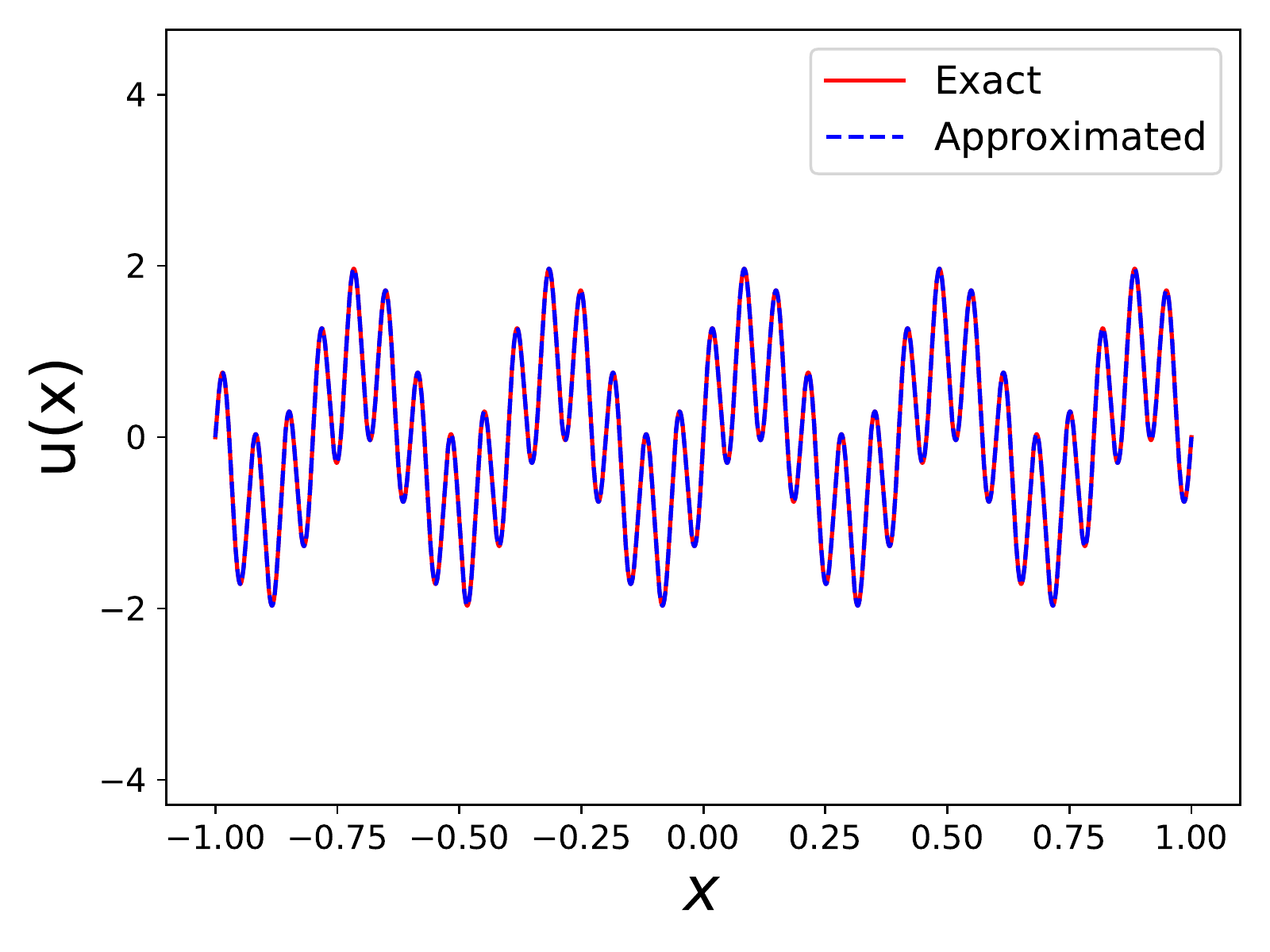}&
    \includegraphics[width=6cm]{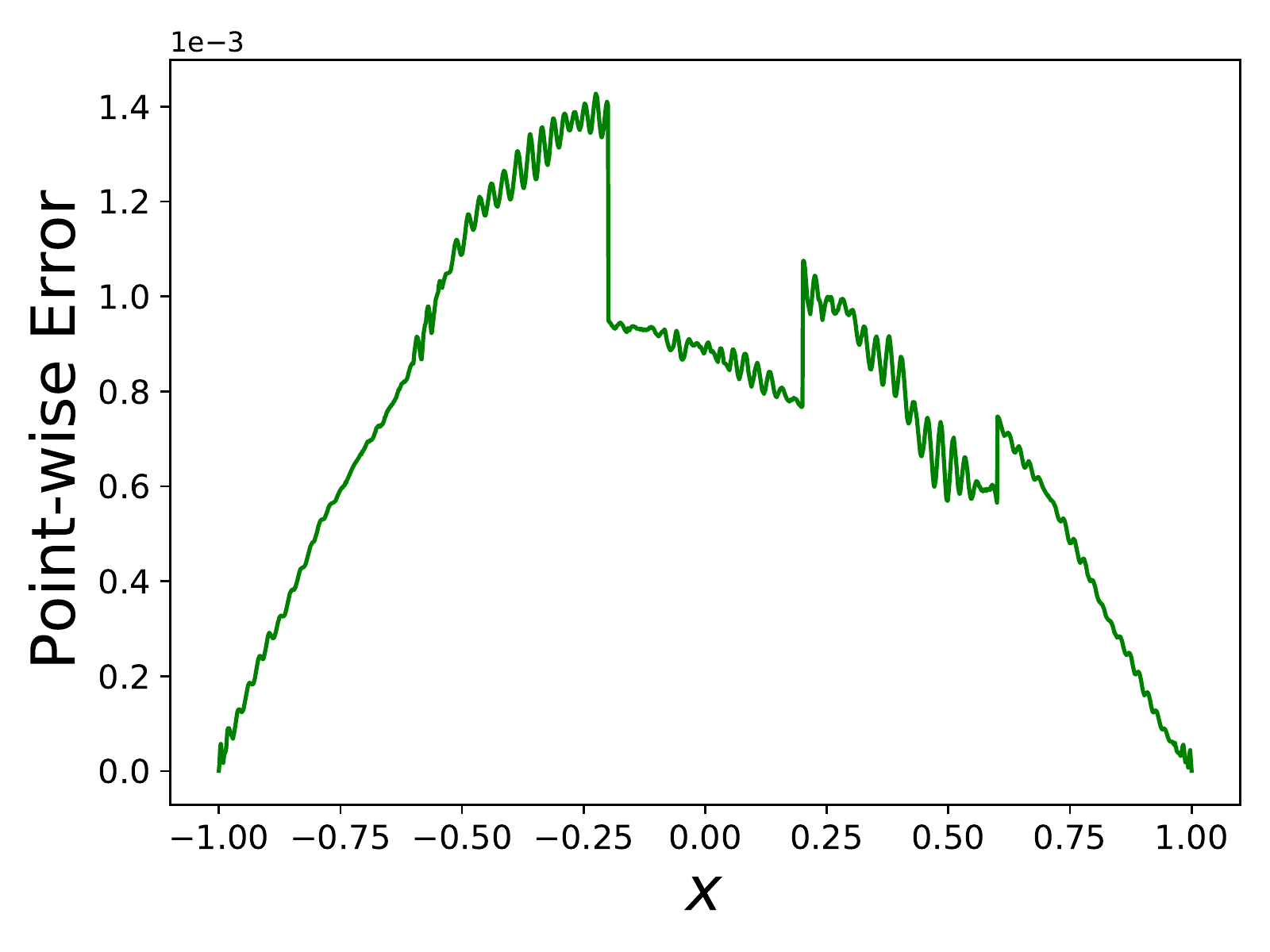}\\
    (e) F-D3M & (f) F-D3M point-wise error
     \end{tabular}}
     \caption{Approximation solutions and point-wise errors,  one-dimensional diffusion problem.}
    \label{fig_1d_poisson}
\end{figure}

To assess the accuracy of neural network approximations, 
2000 test points are uniformly sampled with range $[-1,1]$,
and values of the exact solution and neural network approximations are compared for these test points.
Here, the absolute point-wise errors and the relative errors (see \eqref{relative_l2_error}) are computed,
and all results are repeated five times for different random seeds to obtain average errors.  
Figure~\ref{fig_1d_poisson} shows the results of F-D3M, Global-MFFNet and Global-FCN.
Figure~\ref{fig_1d_poisson}(a), Figure~\ref{fig_1d_poisson}(c) and Figure~\ref{fig_1d_poisson}(e) 
show the solutions of the three methods, where it can be seen that the Global-FCN solution is inaccurate, 
and our F-D3M is very accurate---the F-D3M approximation and the exact solution are visually indistinguishable. 
The absolute point-wise errors of these methods are shown in Figure~\ref{fig_1d_poisson}(b), Figure~\ref{fig_1d_poisson}(d) and Figure~\ref{fig_1d_poisson}(f)
and the relative errors  \eqref{relative_l2_error} of Global-FCN, Global-MFFNet and F-D3M
are $1.86\mathrm{e}+00$,  $5.18\mathrm{e}-01$ and $3.79\mathrm{e}-04$ respectively.
It is clear that Global-FCN has a very large error, and it is consistent with the results in \cite{rahaman2019spectral,ronen2019convergence,xu2020frequency,luo2019theory}, where it is shown that FCN has the so-called ``spectral bias" pathology for high frequency modes. 
While Global-MFFNet has smaller than Global-FCN, the relative error of our F-D3M is three orders 
of magnitude smaller than that of Global-MFFNet. 

\begin{figure}[!ht]
    \centerline{
    \begin{tabular}{cc}
    \includegraphics[width=8cm]{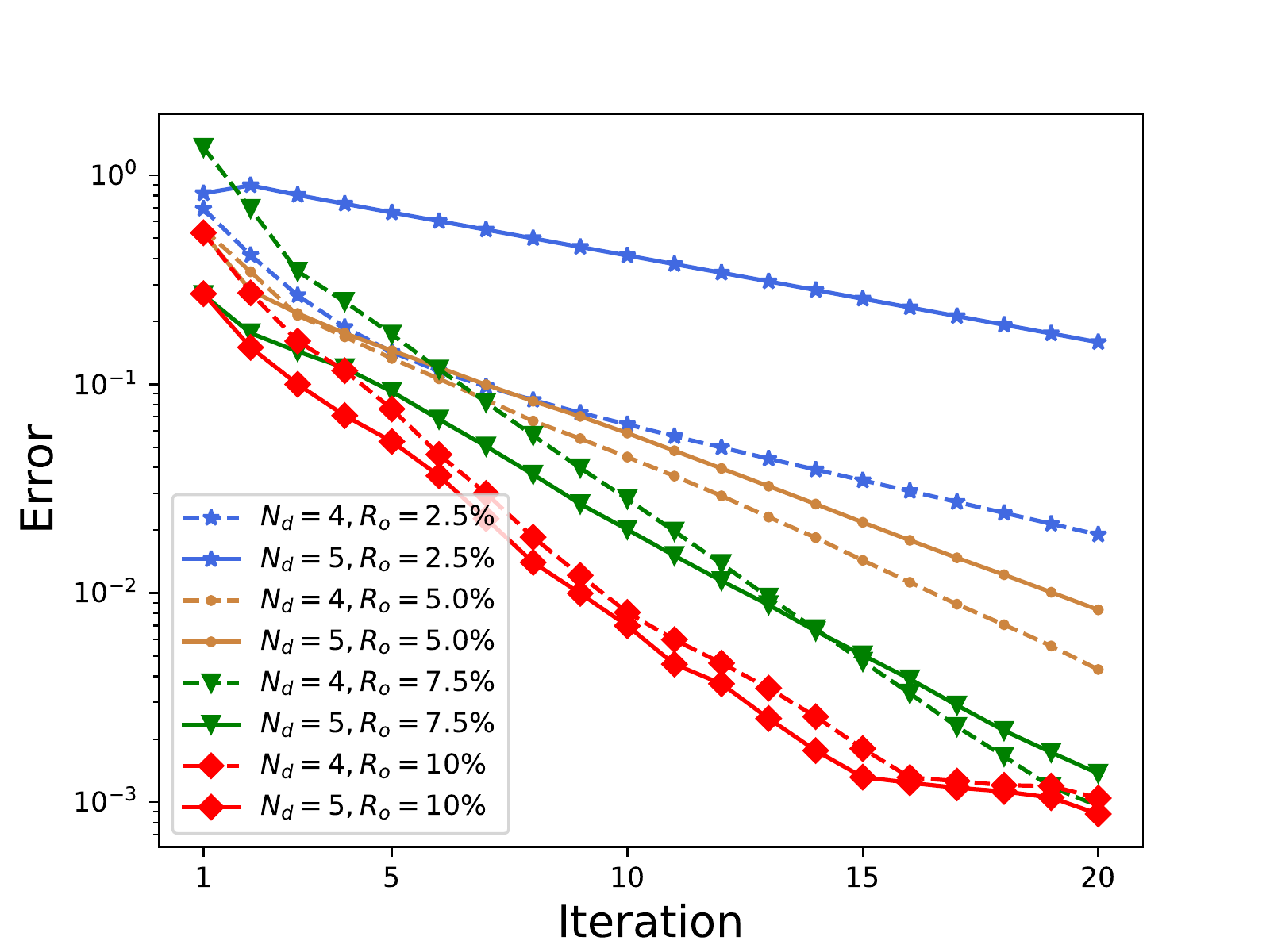}
    \end{tabular}}
    \caption{Errors with respect to domain decomposition iterations, one-dimensional diffusion problem.}
    \label{fig_dd_1d_poisson}
\end{figure}

Effects of different subdomain numbers and different overlapping region sizes are tested next.
Two cases for subdomain numbers are considered---one for $N_d=4$, and the other for $N_d=5$. For each case, the subdomains have the same size and are defined by \eqref{x1}.
Defining the overlapping region proportion as $R_o:=\frac{w}{2}$  where $2$ is the width of the global domain, four different values of $R_o$ are considered, which include $R_0=2.5\%$, $R_0=5\%$
and $R_0=7.5\%$ and $R_0=10\%$. Figure~\ref{fig_dd_1d_poisson} shows the errors with respect to domain decomposition iterations, i.e., $\epsilon(u^{k}(\bx;\btheta),u(\bx))$ (see \eqref{relative_l2_error}),
where $u^{k}(\bx;\btheta)$ is the approximation obtained at iteration step $k$ (line  23 of Algorithm \ref{alg_F_D3M}) and $u(\bx)$ is the solution \eqref{test1-sol}. 
It can be seen that the errors reduce as the number of iterations increases. 
For the case $N_d=4$ (or for $N_d=5$), 
it is clear that large overlapping region proportions (larger values for $R_o$) lead to 
small errors, which is expected for overlapping domain decomposition methods. 
For each value of the overlapping region proportion, errors for the two cases of subdomain numbers ($N_d=4$ and $N_d=5$) are similar, expect for the situation that the overlapping  proportion is very small ($R_o=2.5\%$) 
where the errors for both subdomain numbers are large.
So, compared to the number of subdomains, the size of overlapping regions plays a more important role 
for F-D3M, which is also expected for an overlapping domain decomposition based strategy.
\subsection{Two-dimensional diffusion problem}
For this test problem, 
the following two-dimensional Poisson equation is considered,
\begin{equation}
    \begin{aligned}
        \Delta u(x_1,x_2)&=f(x_1,x_2)\quad \text{in}\ \Omega, \\
        u(x_1,x_2)&=g(x_1,x_2)\quad  \text{on} \ \pOmega,
    \end{aligned}
    \label{test2-pde}
\end{equation}
where the spatial domain is $\Omega=[-1,1]^2$ and the exact solution considered is  
\begin{equation}
    \begin{aligned}
        u(x_1,x_2)=\sin(5\pi x_1)\cos(3\pi x_2)+\sin(10\pi x_1)\cos(2\pi x_2)+\sin(20\pi x_1)\cos(\pi x_2).\\
    \end{aligned}
    \label{test2-sol}
\end{equation}
In \eqref{test2-pde}, $f$ and $g$ are specified by \eqref{test2-sol}. It is clear that the solution \eqref{test2-sol} contains multiple high frequency modes.

\begin{figure}[!ht]
    \centerline{
    \begin{tabular}{cc}
    \includegraphics[width=5cm]{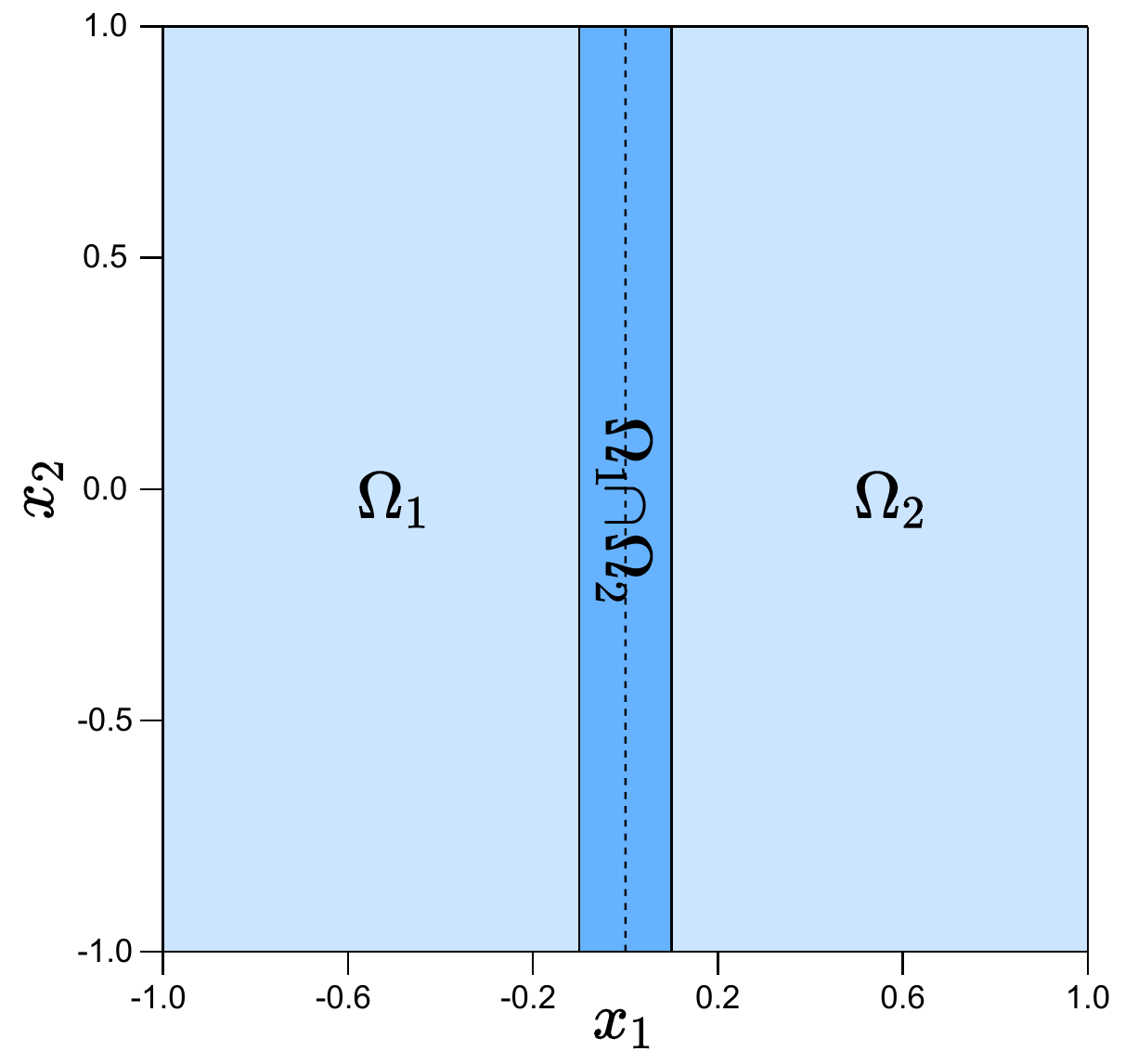}&
        \includegraphics[width=5cm]{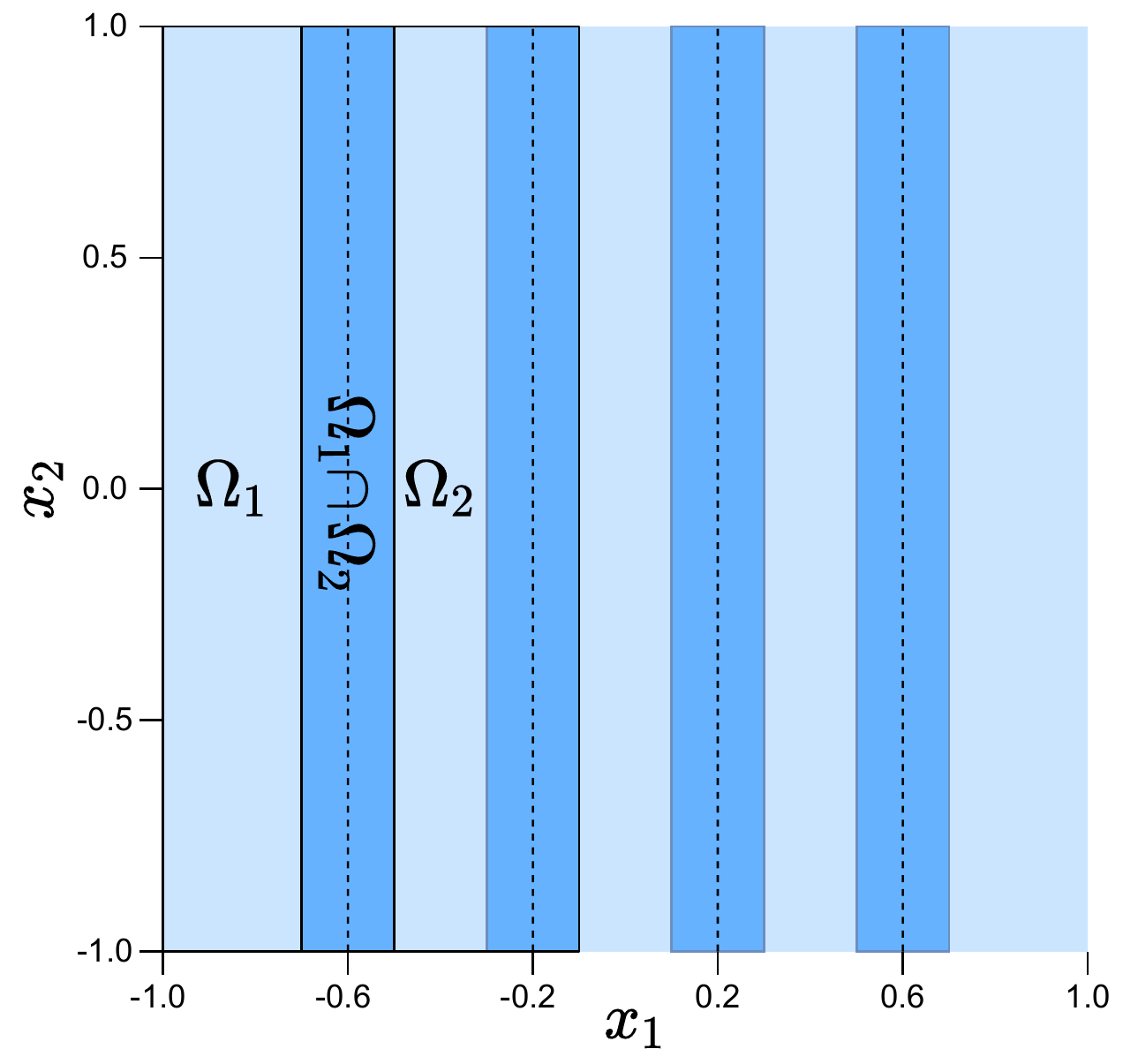}   \\
    (a)  $N_d=2$ subdomains & (b) $N_d=5$ subdomains\\
    \end{tabular}}
    \caption{Overlapping domain decomposition, two-dimensional diffusion problem.}
    \label{fig_dd_definition}
\end{figure}

The spatial domain $\Omega$ is equally divided into $N_d$ overlapping subdomains, and 
two cases with $N_d=2$ and $N_d=5$ are considered, which 
are shown in Figure~\ref{fig_dd_definition}(a) and Figure~\ref{fig_dd_definition}(b) respectively. 
The width of the overlapping region is set to $w=0.2$ for the two cases,  
each subdomain is defined as $\{\Omega_i\}_{i=1}^{N_d}=\{[x_{left}^{(i)},x_{right}^{(i)}]\times[-1,1]\}_{i=1}^{N_d}$ where $x_{left}^{(i)}$ and $x_{right}^{(i)}$ are defined in \eqref{x1}.  
For these regular rectangular subdomains,
 the interface functions 
only depend on $x_2$ (where $\bx=[x_1,x_2]^{\text{T}}$), 
and we denote them by $\phi_{i,j}^{k}(x_2)$.
In addition,  
$\Phi_i^{k}(\bx)$ and $D_i(\bx)$ can be written analytically in this situation. 
For example, considering the subdomain $\Omega_1$ (in the left of $\Omega$),
we denote it as $\Omega_1=[a,b]\times[c,d]$ for simplicity.
Its has three exterior boundaries, and the corresponding boundary functions are
denoted by $g_u(x_1)=g(x_1,x_2)|_{x_2=c}$, $g_t(x_1)=g(x_1,x_2)|_{x_2=d}$ and $g_l(x_1)=g(x_1,x_2)|_{x_2=c}$ respectively. The interface function of $\Omega_1$ is denoted by $\phi_{1,2}^k(x_2)$.  
Then the boundary conditions of the local problem posed on  $\Omega_1$ are
\begin{eqnarray*}
        u_1^{k}(x_1,x_2)&=&g_u(x_1)\quad \text{on}\, [a,b]\times\{c\},\\
        u_1^{k}(x_1,x_2)&=&g_t(x_1)\quad \text{on}\, [a,b]\times\{d\},\\
        u_1^{k}(x_1,x_2)&=&g_l(x_2)\quad \text{on}\, \{a\}\times[c,d],\\
        u_1^{k}(x_1,x_2)&=&\phi_{1,2}^{k}(x_2)\quad \text{on}\, \{b\}\times[c,d],
    \label{bound_conditions}
\end{eqnarray*}
where $\{\bx|x_1=b,\,c\leq x_2 \leq d\}$ is an interface.
Then $\Phi_1^{k}(\bx)$ and $D_1(\bx)$ can be written as 
\begin{equation*}
    \begin{aligned}
        \Phi_1^{k}(\bx)=&\frac{1}{\omega}(b-x_1)g_l(x_2)+\frac{1}{\omega}(x_1-a)\phi_{1,2}^{k}(x_2)\\
        &+\frac{1}{h}(d-x_2)\left\{g_u(x_1)-\left[\frac{1}{\omega}(b-x_1)g_u(a)+\frac{1}{\omega}(x_1-a)g_u(b)\right]\right\}\\
        &+\frac{1}{h}(x_2-c)\left\{g_t(x_1)-\left[\frac{1}{\omega}(b-x_1)g_t(a)+\frac{1}{\omega}(x_1-a)g_t(b)\right]\right\},
    \end{aligned}
\end{equation*}
and
\begin{equation*}
    D_1(\bx)=(x_1-a)(x_1-b)(x_2-c)(x_2-d),
\end{equation*}
where $\omega=b-a$ and $h=d-c$. 
For other subdomains $\Omega_i$ for $i=2,\ldots,N_d$, $\Phi_i^{k}(\bx)$ and $D_i(\bx)$ can be obtained similarly.

In our F-D3M algorithm (Algorithm \ref{alg_F_D3M}), a MFFNet with four RFF subnetworks is applied for each subdomain for this test problem. The standard deviation $\sigma$ is set to $1,5,10,20$ for the four subnetworks respectively, and each subnetwork has two hidden layers with fourty neurons in each layer. We also consider the aforementioned Global-FCN and Global-MFFNet methods for comparison and the weight coefficient in \eqref{em_loss_func} is set to $\lambda=100$ for both methods. The setting of MFFNet in Global-MFFNet is the same as that for F-D3M, while Global-FCN has two hidden layers with one hundred sixty neurons in each layer. Networks in the three methods then have the same number of neurons. For all networks, the Kaiming scheme \cite{he2015delving} is used for initialization and the hyperbolic tangent function is employed as the activation function. Moreover, all networks are trained by the Adam \cite{kingma2014adam} optimizer with default settings. An exponential learning rate with the initial learning rate of 0.01 and a decay rate of 0.9 every 500 training epochs is employed, while the learning rate is reinitialized for each domain decomposition iteration in F-D3M.  
As the local solutions in F-D3M gradually get close to the global solution as the domain decomposition iteration step
increases, we gradually increase the number of  training epochs---at the initial step, the number of epochs
is set to $2500$, and it increases by $500$ for each step until 
the number of domain decomposition iterations reaches 15.   
For each training epoch, $|\mS_i|=1000$ collocation points are randomly sampled using the uniform distribution with range $\Omega_i$ for $i=1,\ldots,N_d$. 
To train Global-FCN and Global-MFFNet, $90000$ epochs are used and 
the number of collocation points for the interior domain and the boundary 
are set to $N_r=5000$ and $N_b=200\times4$ respectively.

\begin{figure}[!ht]
    \centerline{
    \begin{tabular}{c}
    \includegraphics[width=6cm]{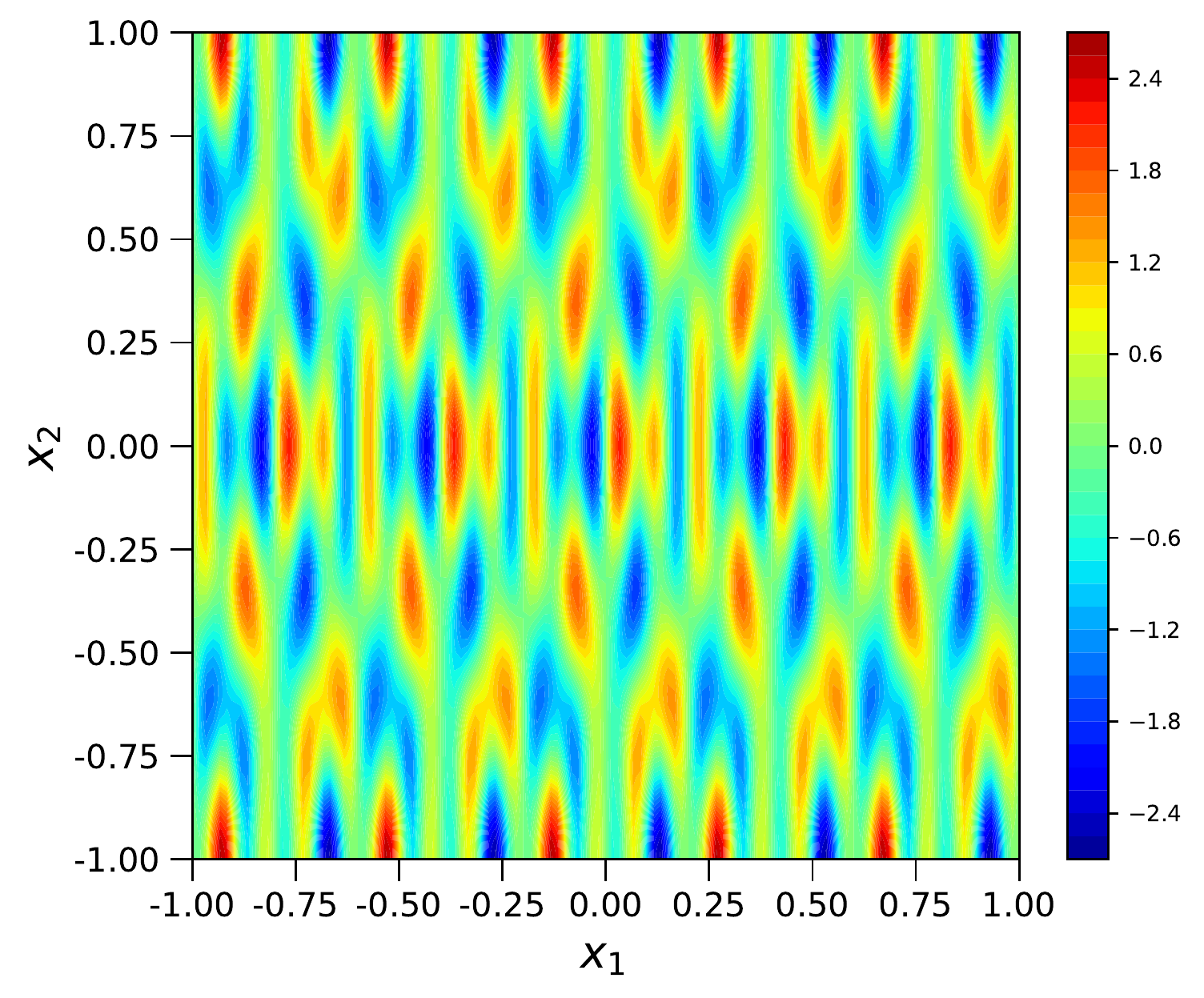}\\
    \end{tabular}}
    \caption{Exact solution, two-dimensional diffusion problem.}
    \label{fig_2d_poisson_exact}
\end{figure}

\begin{figure}[!ht]
    \centerline{
    \begin{tabular}{ccc}
    \includegraphics[width=5cm]{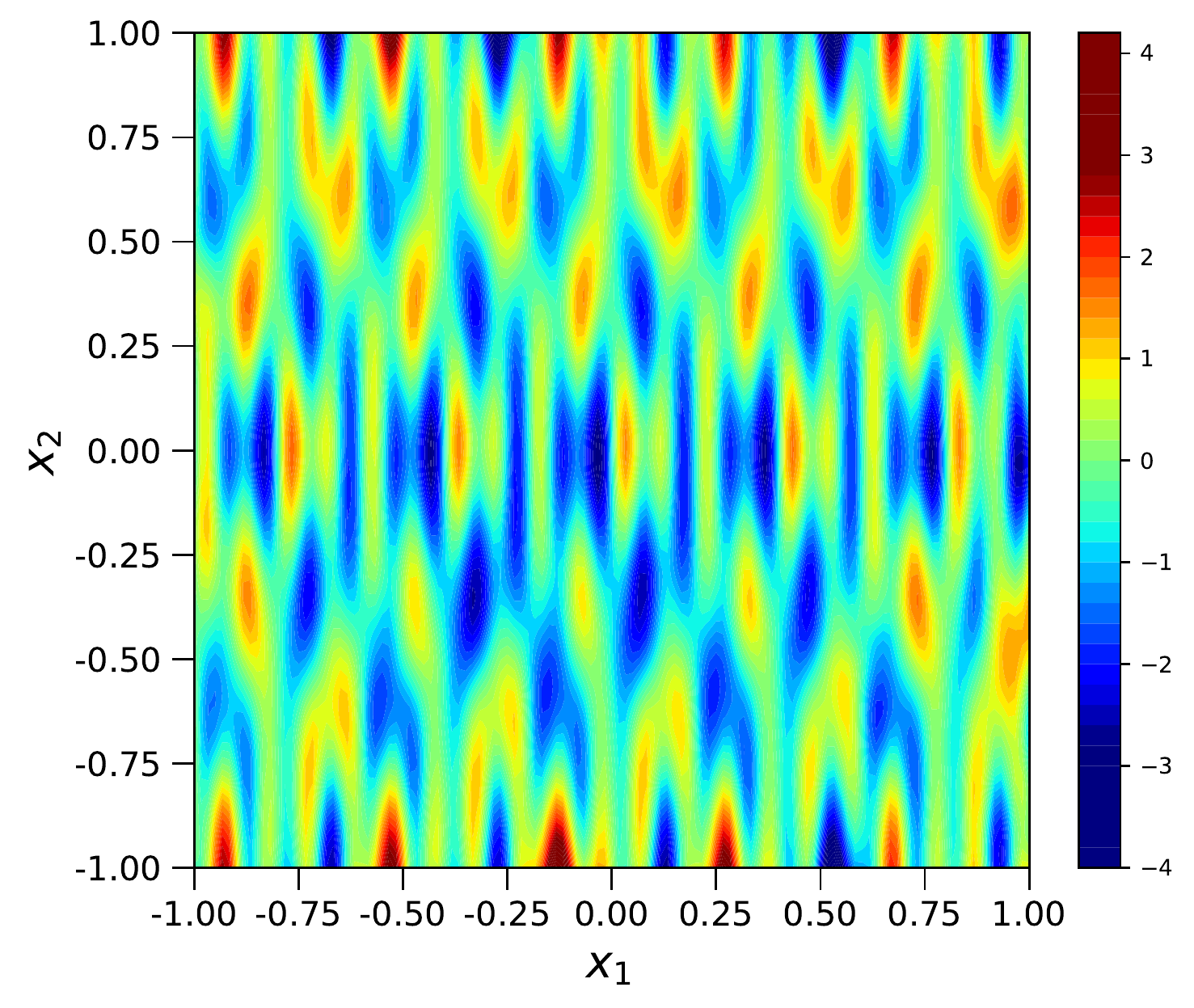}&\quad\quad\quad
    \includegraphics[width=5cm]{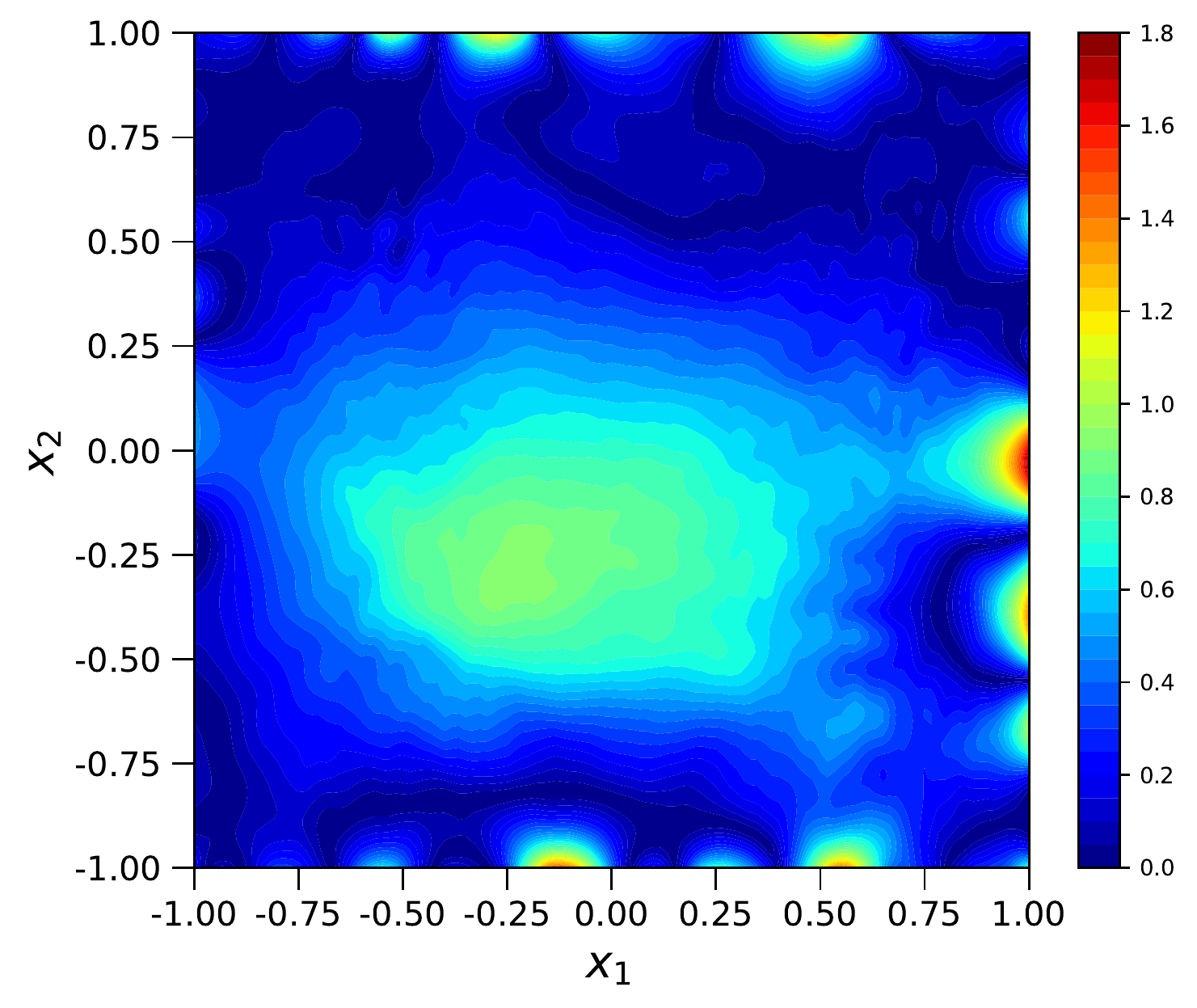}\\
    (a) Global-FCN appoximation &\quad\quad\quad (b) Global-FCN point-wise error\\
    \includegraphics[width=5cm]{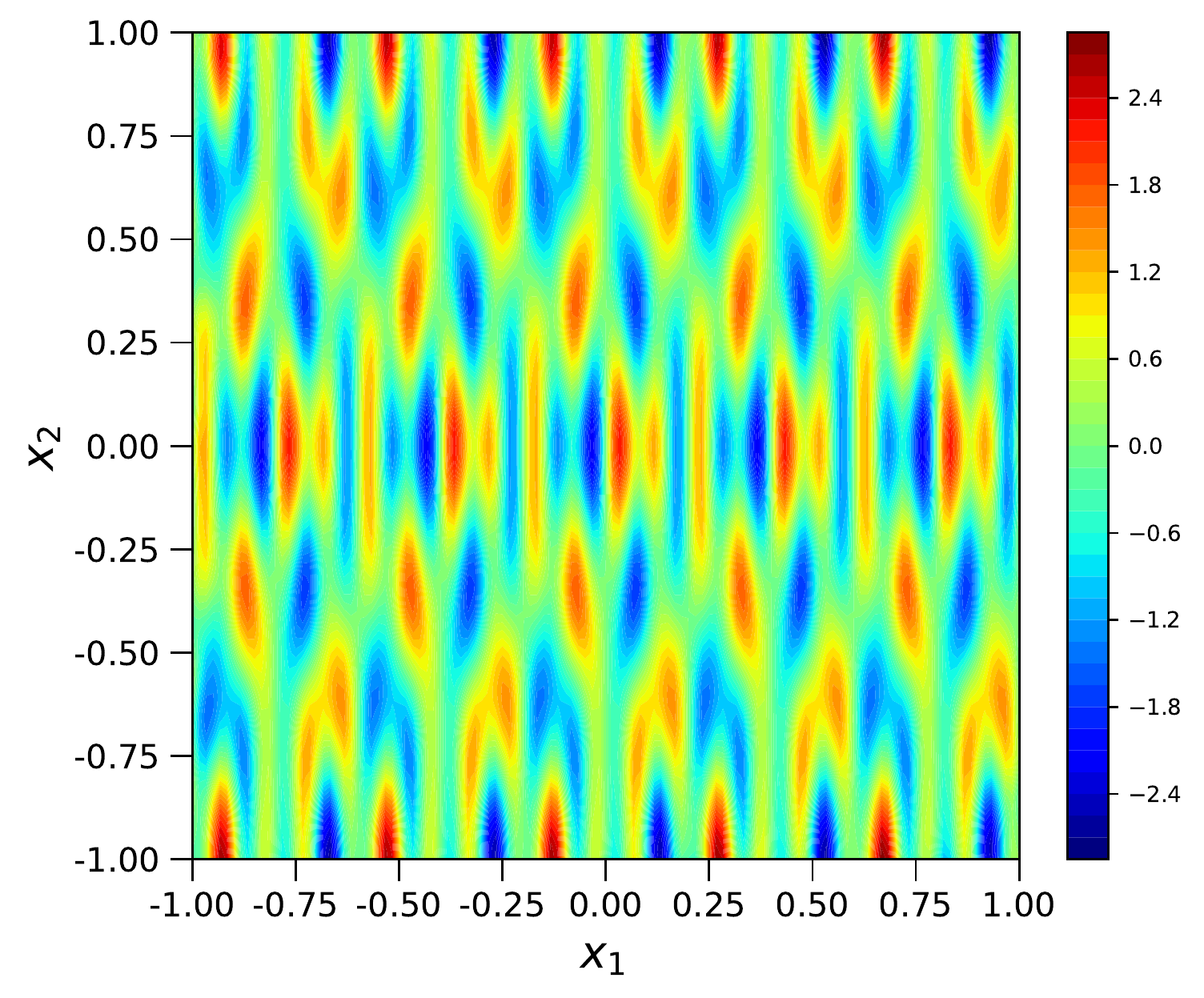}&\quad\quad\quad
    \includegraphics[width=5cm]{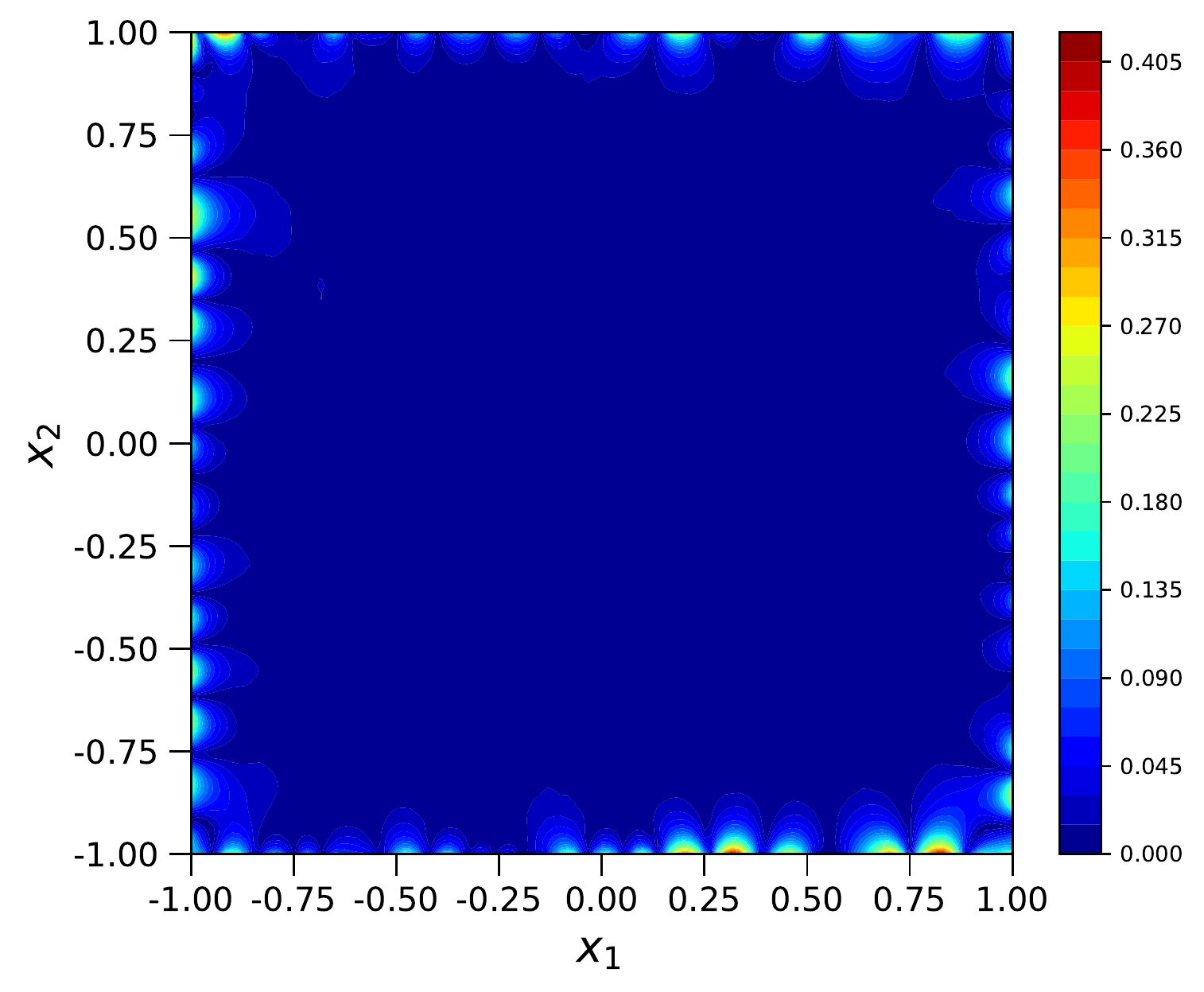}\\
    (c) Global-MFFNet appoximation &\quad\quad\quad (d) Global-MFFNet point-wise error\\
    \includegraphics[width=5cm]{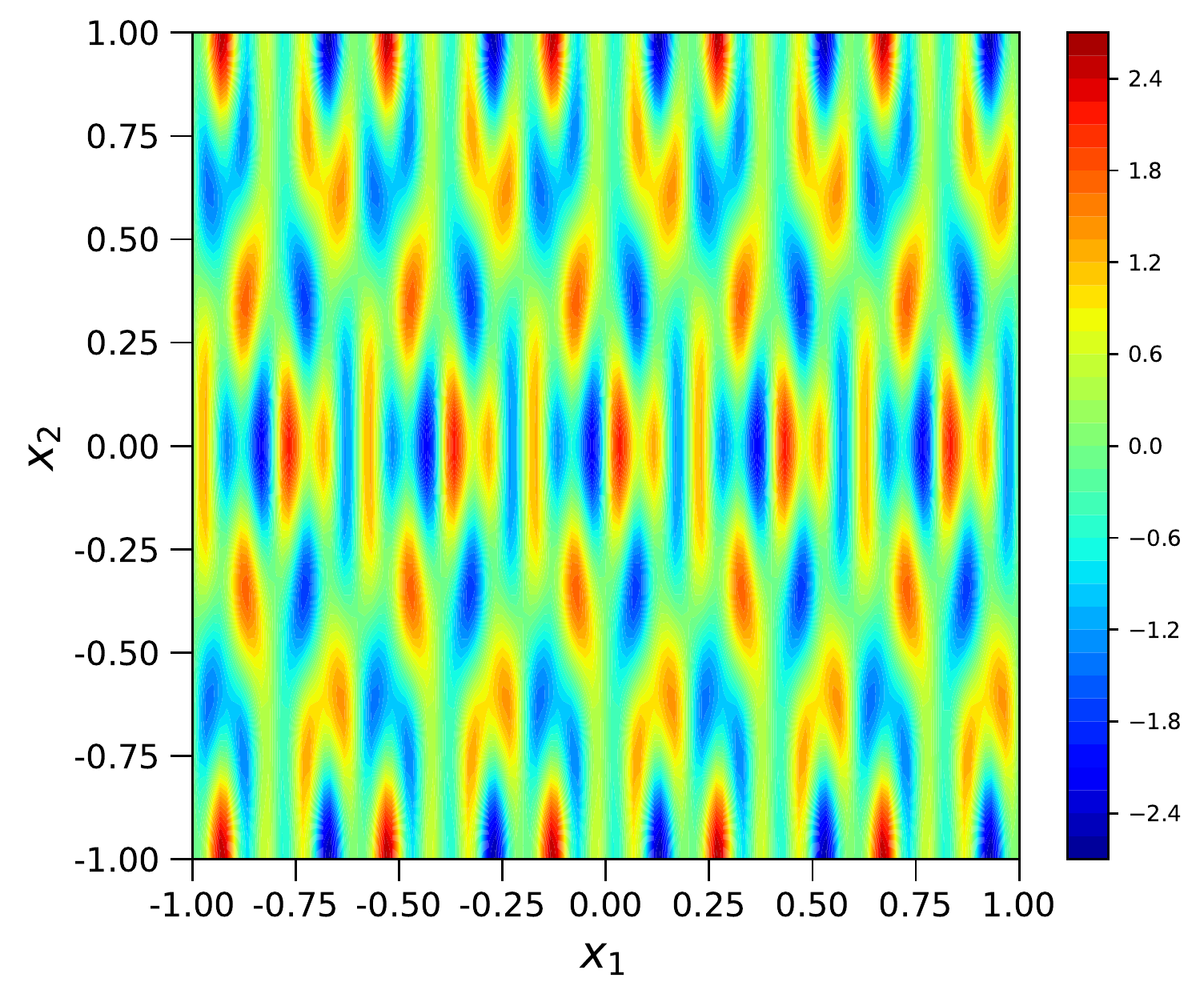}&\quad\quad\quad
    \includegraphics[width=5cm]{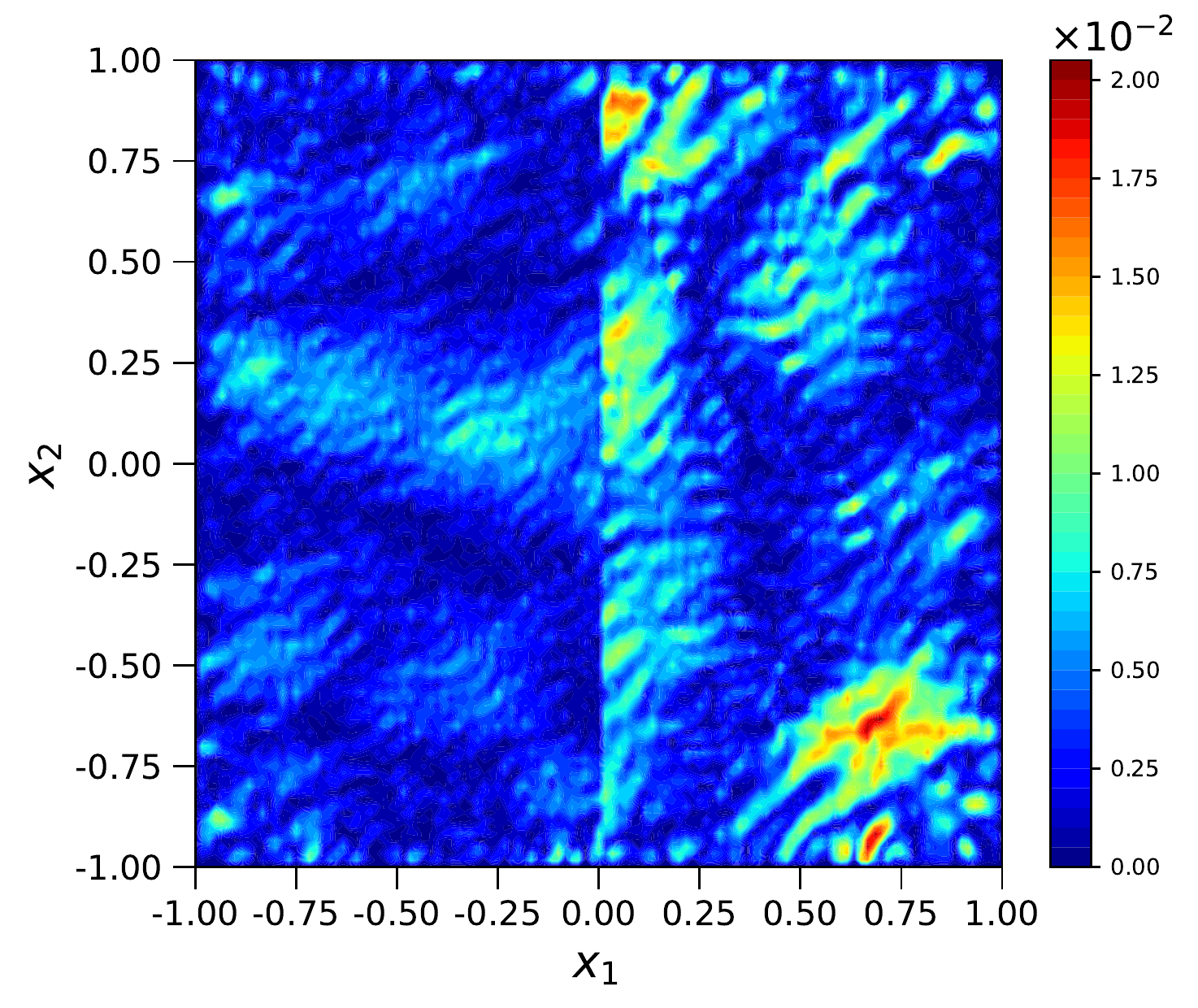}\\
    (e) F-D3M ($N_d=2$) appoximation &\quad\quad\quad (f) F-D3M ($N_d=2$) point-wise error\\
    \includegraphics[width=5cm]{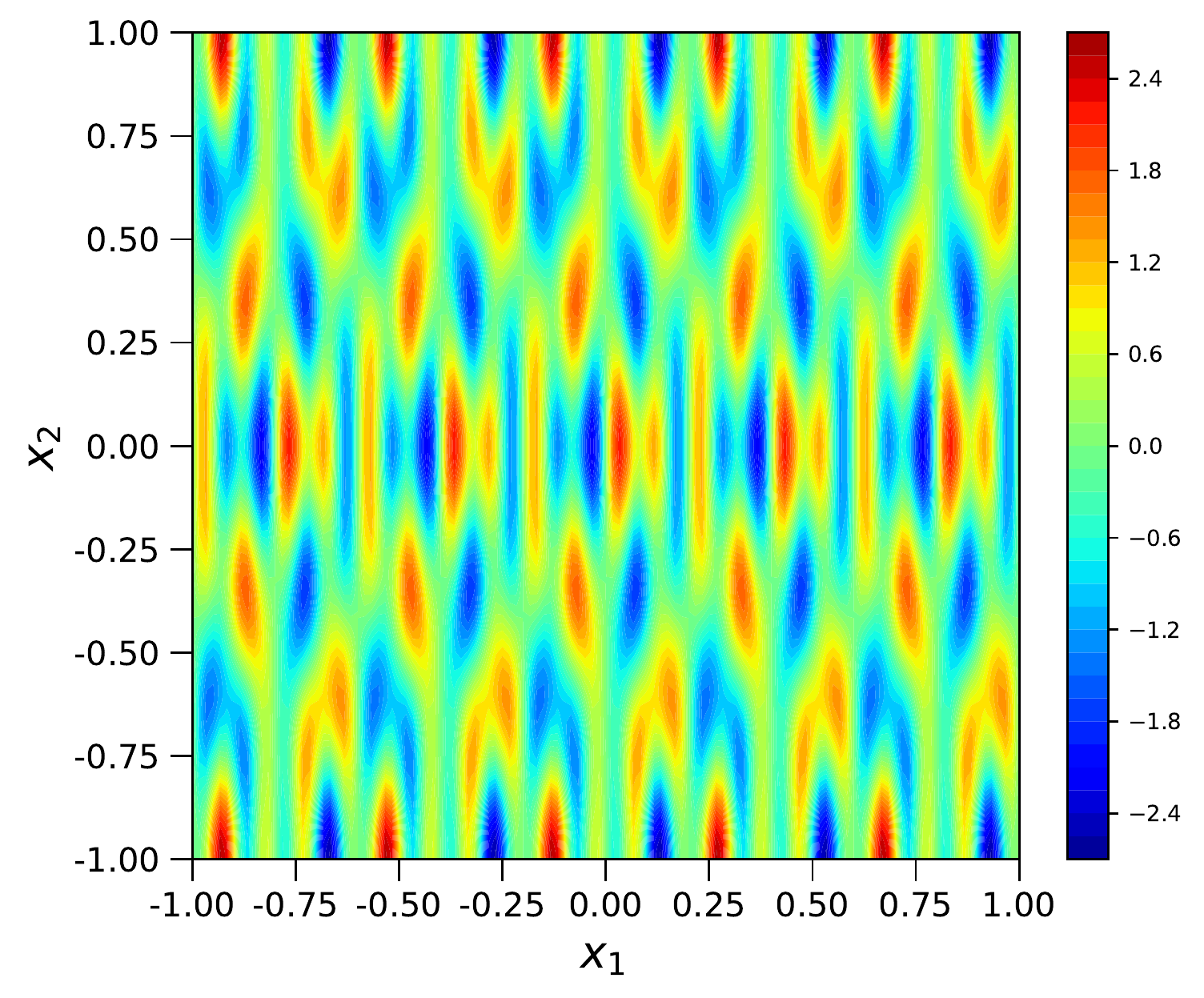}&\quad\quad\quad
    \includegraphics[width=5cm]{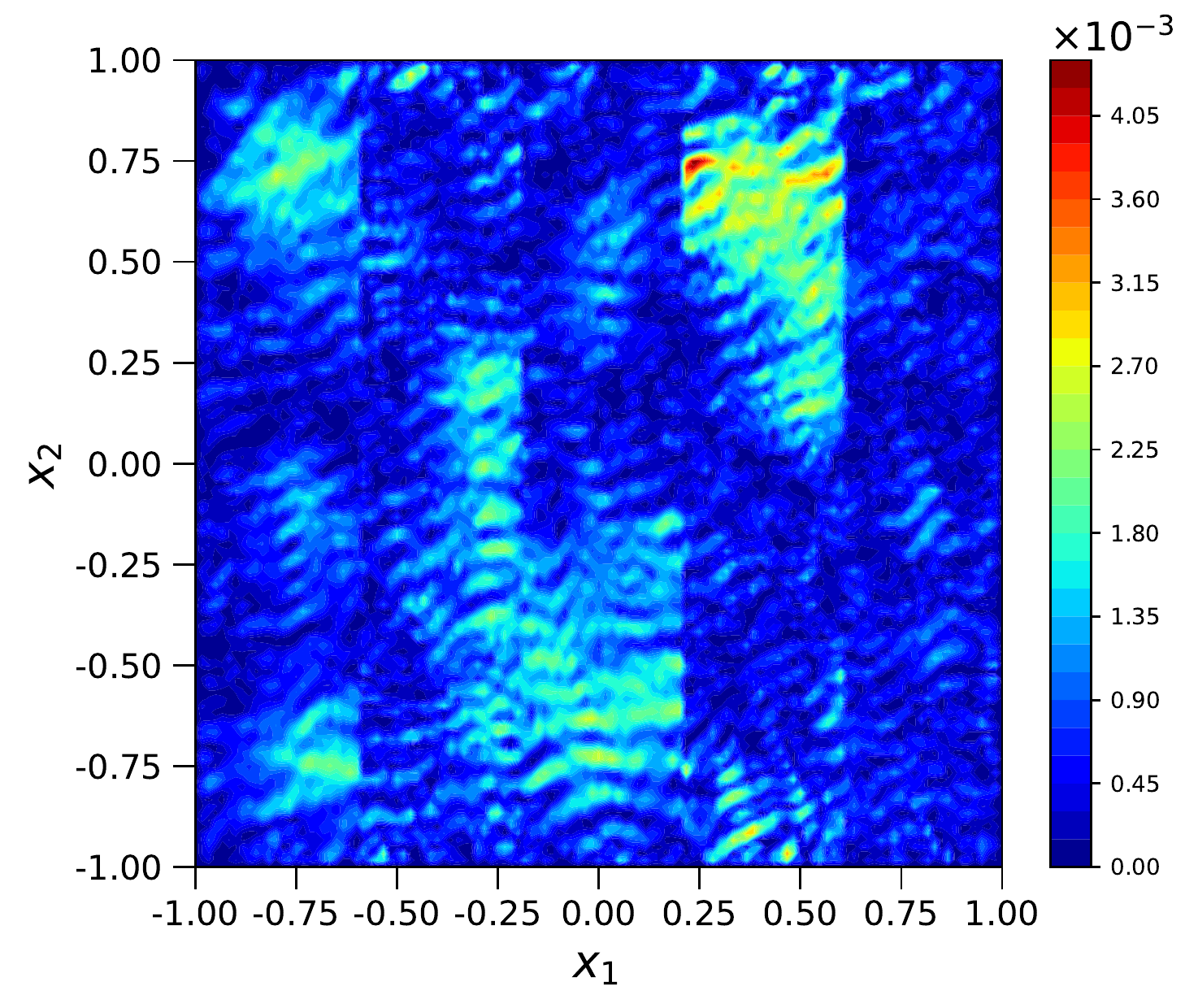}\\
    (g) F-D3M ($N_d=5$) appoximation &\quad\quad\quad (h) F-D3M ($N_d=5$) point-wise error\\
    \end{tabular}}
    \caption{Approximation solutions and point-wise errors,  two-dimensional diffusion problem.}
    \label{fig_2d_poisson}
\end{figure}

To evaluate the accuracy of neural network approximations, uniform $121\times 121$ grids for $\Omega$ are used,  and values of the exact solution and neural network approximations are compared for these grid points. 
The absolute point-wise errors and the relative errors (see \eqref{relative_l2_error}) are computed, and all results are repeated three times for different random seeds to obtain average errors. Figure~\ref{fig_2d_poisson_exact} shows the exact solution \eqref{test2-sol} and Figure~\ref{fig_2d_poisson} shows approximation solutions and 
errors of F-D3M with $N_d=2$, F-D3M with $N_d=5$, Global-MFFNet and Global-FCN. 
From Figure~\ref{fig_2d_poisson}(a), Figure~\ref{fig_2d_poisson}(c), Figure~\ref{fig_2d_poisson}(e) and Figure~\ref{fig_2d_poisson}(g), 
it can be seen that the Global-FCN approximation is clearly different from the exact solution shown in  \eqref{test2-sol},
and the approximations obtained by Global-MFFNet and F-D3M with $N_d=2$ and $N_d=5$ are very similar to 
the exact solution.  
The absolute point-wise errors of these methods are shown in Figure~\ref{fig_2d_poisson}(b), Figure~\ref{fig_2d_poisson}(d), Figure~\ref{fig_2d_poisson}(f), Figure~\ref{fig_2d_poisson}(h) 
and the relative errors \eqref{relative_l2_error} of Global-FCN, Global-MFFNet, F-D3M with $N_d=2$ and F-D3M with $N_d=5$ are $4.97\mathrm{e}-01$, $3.93\mathrm{e}-02$, $5.23\mathrm{e}-03$ and $2.70\mathrm{e}-03$, respectively. It is clear that Global-FCN has the largest error. While Global-MFFNet is more accurate  than Global-FCN, the relative errors of our F-D3M are significantly smaller than that of Global-MFFNet.  
The relative errors  \eqref{relative_l2_error} associated with different network sizes are shown in Table~\ref{table_1d_poisson}, where $l$ and $n$ represent the number of hidden layers and the number of neurons in each layer respectively.  
These results show that  for a given network size,  F-D3M has smaller errors
than Global-FCN and Global-MFFNet. 
\begin{table}[!ht]
    \begin{center}
    \caption{Relative errors of Global-FCN, Global-MFFNet, F-D3M  with $N_d=2$ and $N_d=5$, 
    associated with different network sizes, two-dimensional diffusion problem.}
    \begin{tabular}{llllll}
    \hline
    Network size & Method &  & & \\
    \cline { 2 - 5 }
    ($l\times n$) & Global-FCN & Global-MFFNet & F-D3M ($N_d=2$) & F-D3M ($N_d=5$) \\
    \hline 
    1$\times$160&$1.17\mathrm{e}+00$&$8.84\mathrm{e}-01$&$5.81\mathrm{e}-02$&$1.69\mathrm{e}-02$\\
    2$\times$160&$4.97\mathrm{e}-01$&$3.93\mathrm{e}-02$&$5.23\mathrm{e}-03$&$2.70\mathrm{e}-03$ \\
    2$\times$120&$5.09\mathrm{e}-01$&$2.11\mathrm{e}-01$&$5.75\mathrm{e}-03$&$2.78\mathrm{e}-03$ \\
    3$\times$120&$2.30\mathrm{e}-01$&$3.59\mathrm{e}-02$&$4.17\mathrm{e}-03$&$2.02\mathrm{e}-03$ \\
    \hline
    \end{tabular}
    \label{table_1d_poisson}
    \end{center}
\end{table}

\begin{figure}[!ht]
    \centerline{
    \begin{tabular}{ccc}
    \includegraphics[width=5cm]{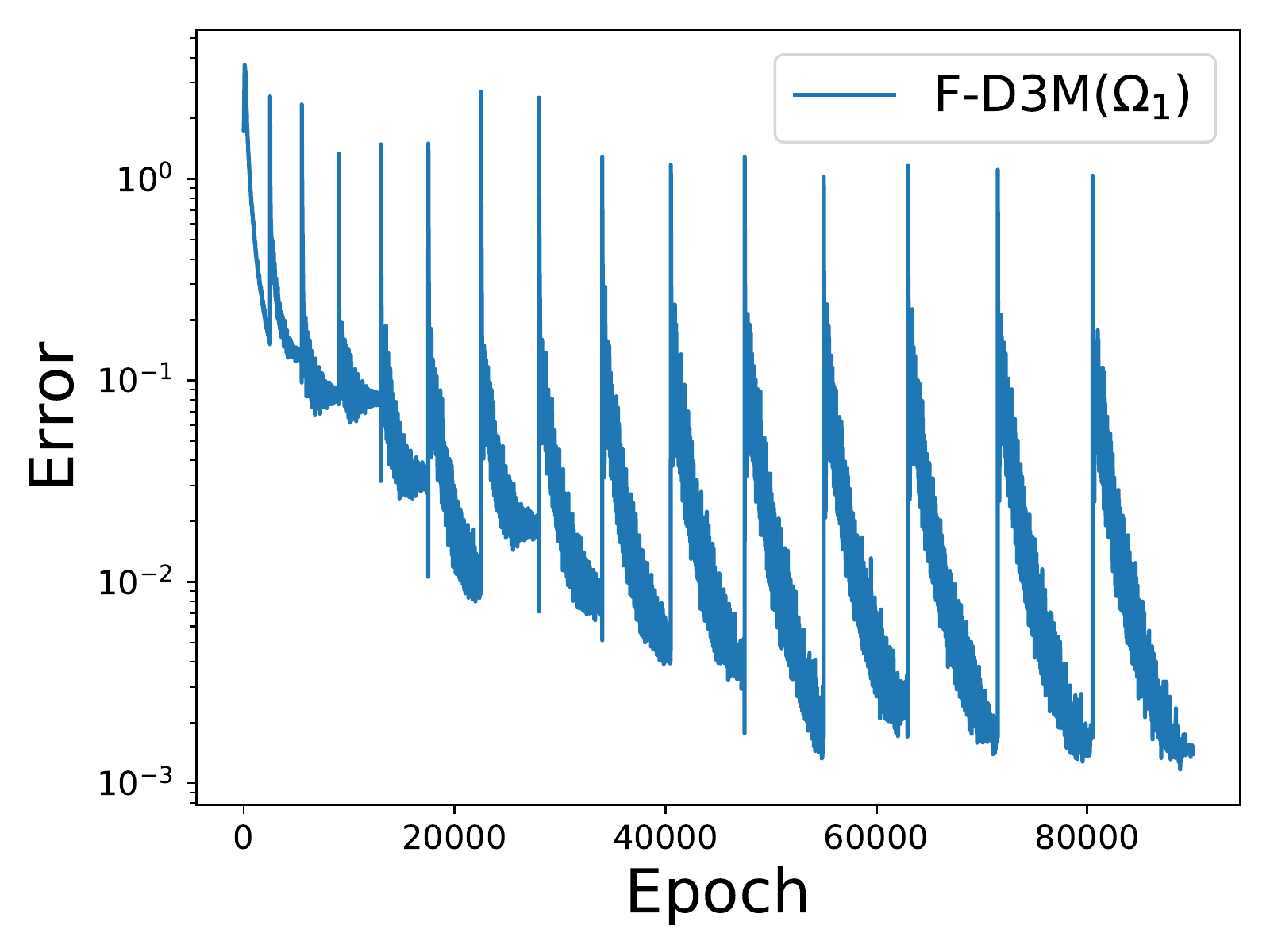}&
    \includegraphics[width=5cm]{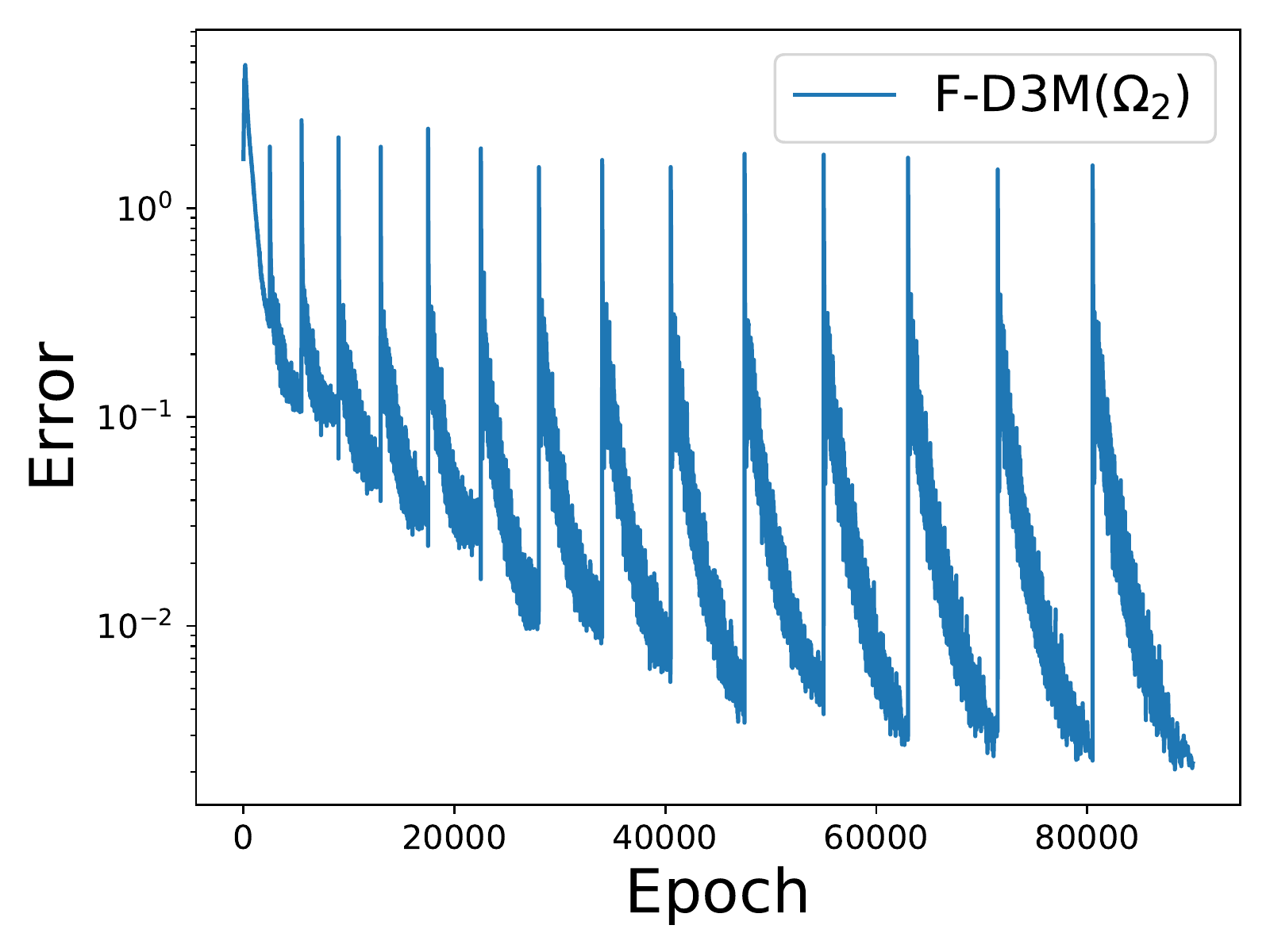}\\
    (a) Errors for $\Omega_1$ &(b) Errors for $\Omega_2$\\
    \includegraphics[width=5cm]{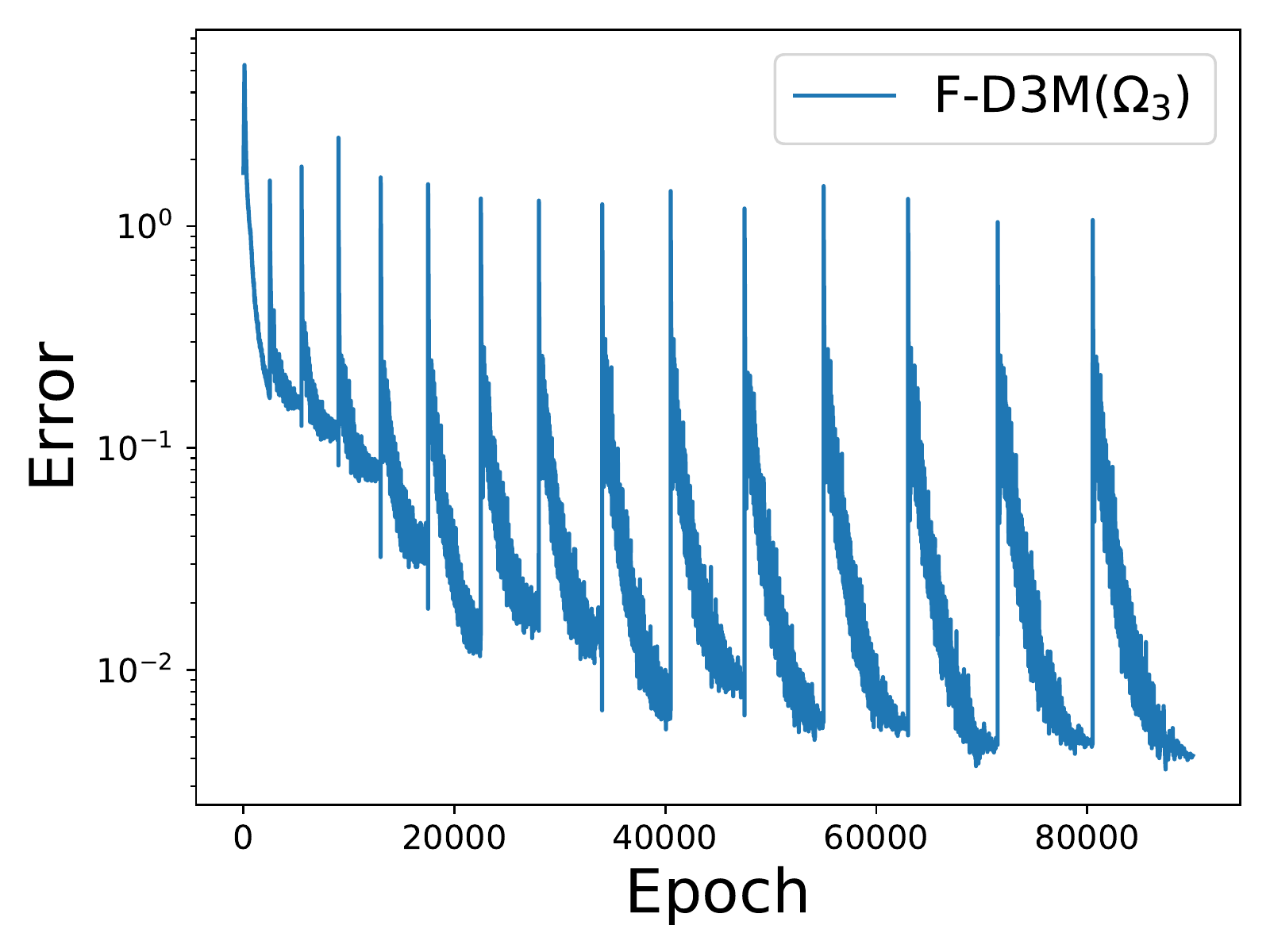}&
    \includegraphics[width=5cm]{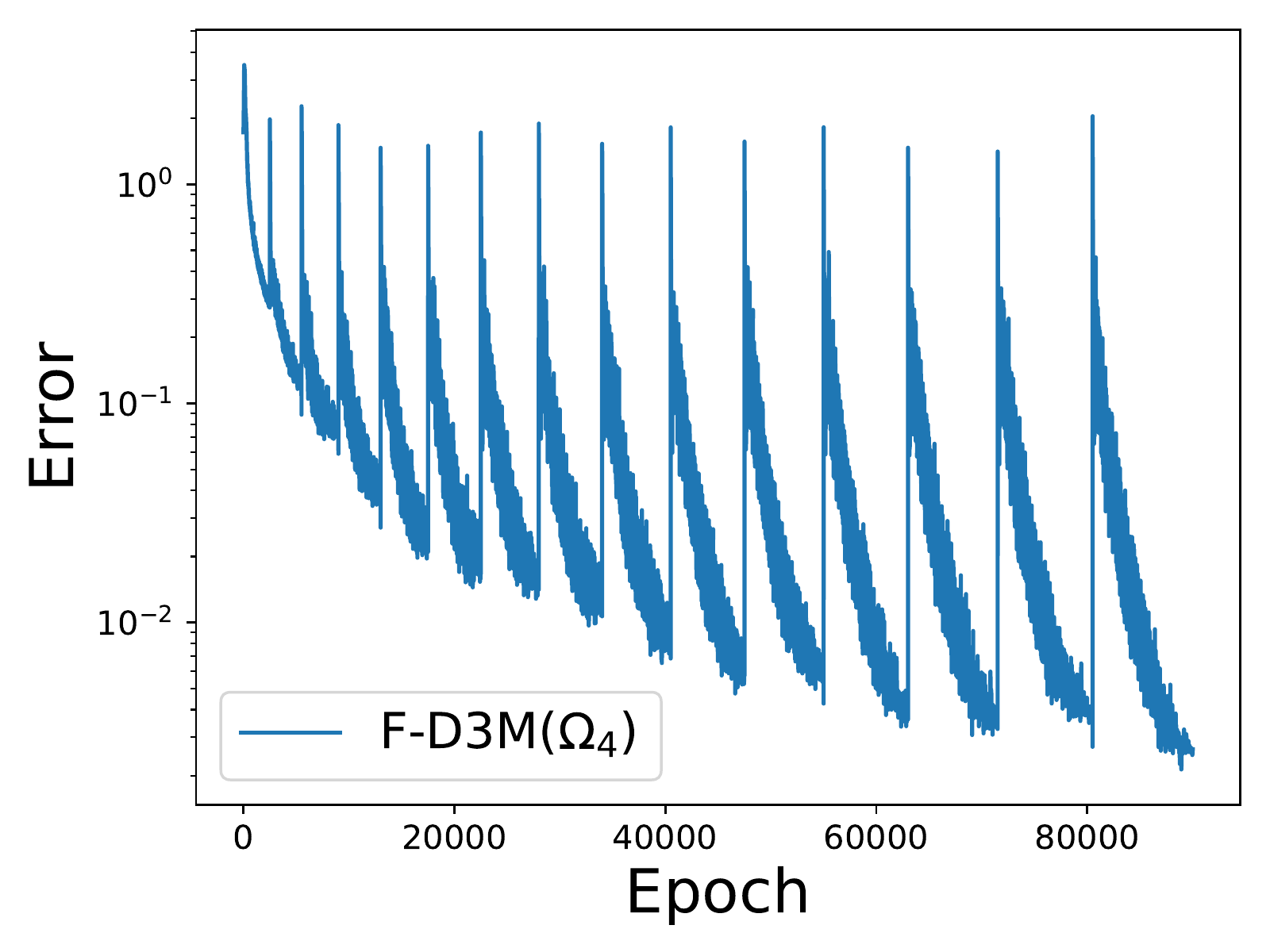}\\
    (c) Errors for $\Omega_3$ &(d) Errors for $\Omega_4$\\
    \includegraphics[width=5cm]{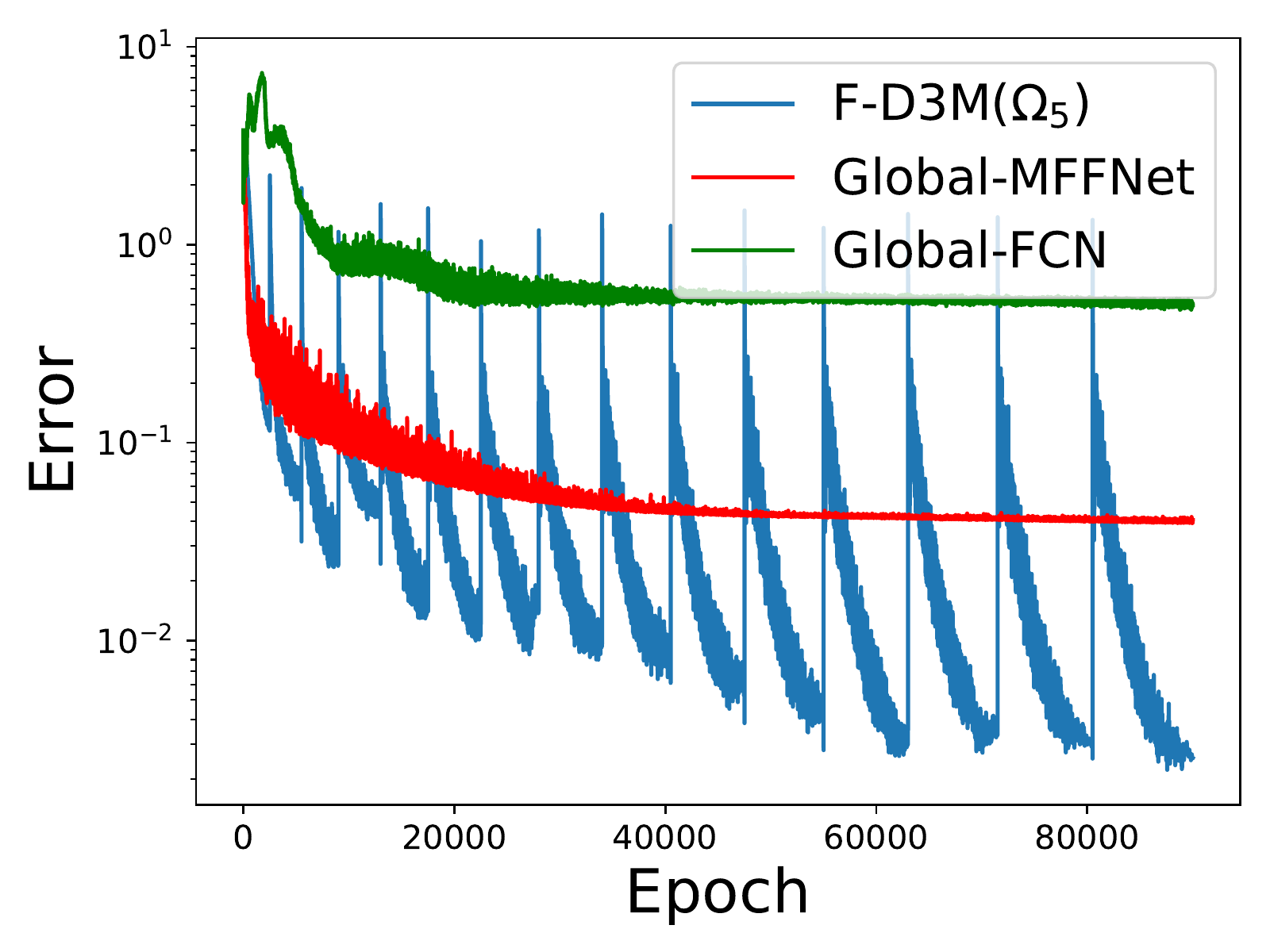}&
    \includegraphics[width=5cm]{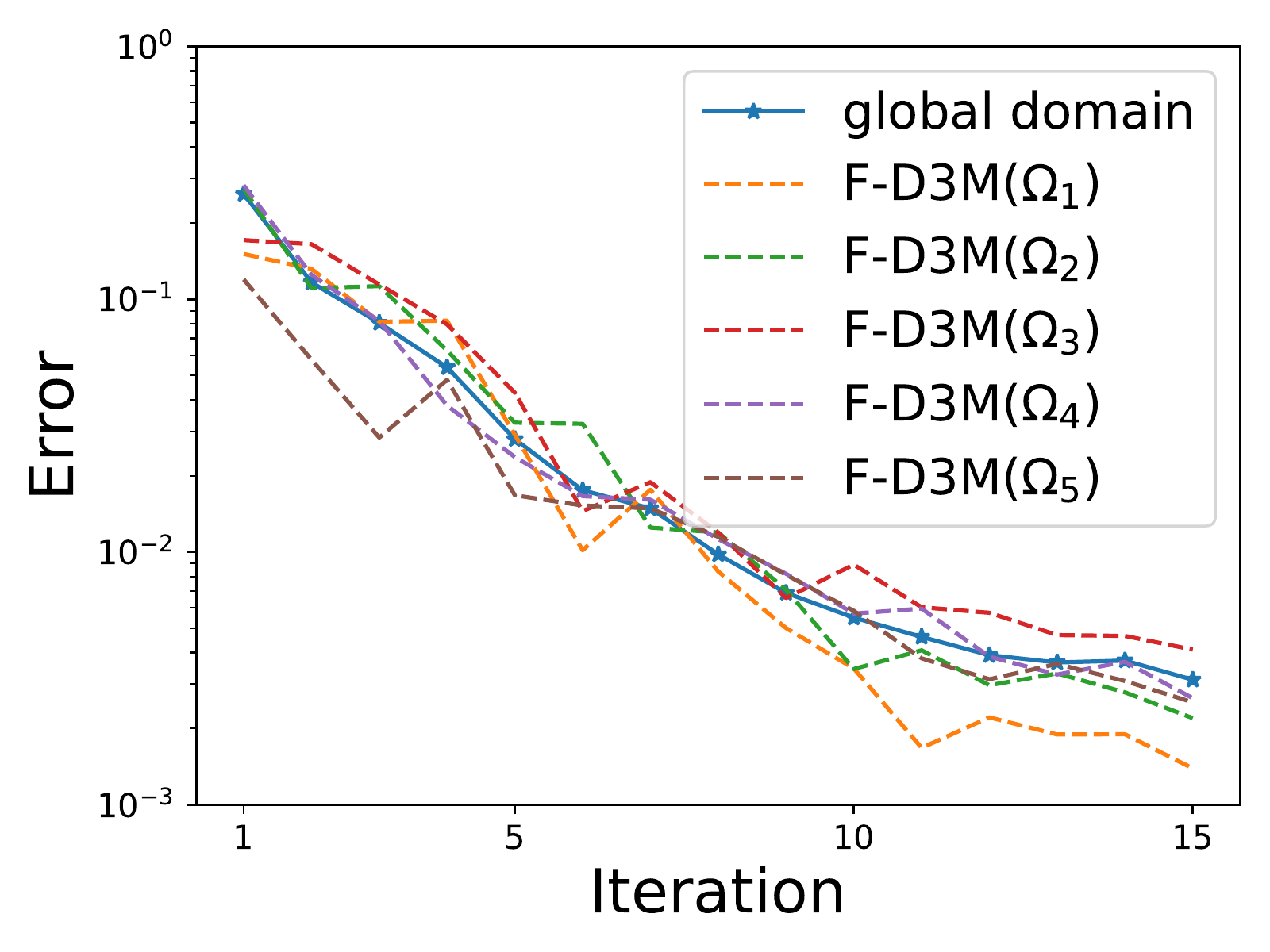}\\
    (e) Errors for $\Omega_5$ &(f) Errors w.r.t. domain decomposition iterations 
    \end{tabular}}
    \caption{Errors w.r.t. training epochs and domain decomposition iterations, two-dimensional diffusion problem.}
    \label{fig_dd_2d_poisson}
\end{figure}

In addition, the relative errors of F-D3M (with $N_d=5$) with respect to training epochs and domain decomposition iteration steps are shown in Figure~\ref{fig_dd_2d_poisson}. 
Figure~\ref{fig_dd_2d_poisson}(a), Figure~\ref{fig_dd_2d_poisson}(b), Figure~\ref{fig_dd_2d_poisson}(c),  
Figure~\ref{fig_dd_2d_poisson}(d) and Figure~\ref{fig_dd_2d_poisson}(e) show the relative errors of F-D3M for subdomains 
$\Omega_1$, $\Omega_2$, $\Omega_3$, $\Omega_4$ and $\Omega_5$ respectively, while 
Figure~\ref{fig_dd_2d_poisson}(e) also shows the relative errors of Global-FCN and Global-MFFNet.
Each sharp jump in Figure~\ref{fig_dd_2d_poisson}(a) to Figure~\ref{fig_dd_2d_poisson}(e) 
is caused by the mismatch between the updated particular solution and MFFNets at first few epochs within each domain decomposition iteration.
During each domain decomposition iteration in F-D3M, the particular solution is updated while MFFNets inherit the parameters from the last iteration directly.
From Figure~\ref{fig_dd_2d_poisson}(e), 
it can be seen that, at each domain decomposition iteration step, as the number of epochs increases,
errors of F-D3M become significantly  smaller than the errors of Global-FCN, Global-MFFNet,
while applying an optimal reinitialization strategy remains an open problem.
Figure~\ref{fig_dd_2d_poisson}(f) shows the relative errors including the error for the global domain $\Omega$ and the errors for each subdomain, reduce as the  domain decomposition iteration  number increases. 
\subsection{Two-dimensional Helmholtz equation}
Due to its wide range of applications in many fields, e.g., electromagnetic radiation, seismology and acoustics,
the  following two-dimensional Helmholtz equation is considered for this test problem, 
\begin{equation}
    \begin{aligned}
        \Delta u(x_1,x_2)+k^2u(x_1,x_2)&=f(x_1,x_2)\quad \text{in}\ \Omega, \\
        u(x_1,x_2)&=g(x_1,x_2)\quad  \text{on} \ \pOmega,
    \end{aligned}
    \label{test3-pde}
\end{equation}
where the spatial domain is $\Omega=[-1,1]^2$ and $k$ is the wavenumber. 
The solution of the Helmholtz equation can have high-frequency oscillations, which 
causes difficulties for applying standard numerical methods, e.g., the finite element methods \cite{elman2014finite}.
The following highly oscillated solution is considered herein,  
\begin{equation}
    \begin{aligned}
        u(x_1,x_2)=e^{\lambda\left(\cos\left(\omega_1 \pi x_1\right)\right)^2}\sin\left(\omega_2\pi x_1\right)\sin\left(\omega_3\pi x_2\right),\\
    \end{aligned}
    \label{test3-sol}
\end{equation}
where we set $k=16  \pi$, $\lambda=\frac{3}{4}$, $\omega_1=2$, $\omega_2=16$ and $\omega_3=1$.
The functions $f(x_1,x_2)$ and  $g(x_1,x_2)$ are specified by \eqref{test3-sol}. 

\begin{figure}[!ht]
    \centerline{
    \begin{tabular}{cc}
    \includegraphics[width=6cm]{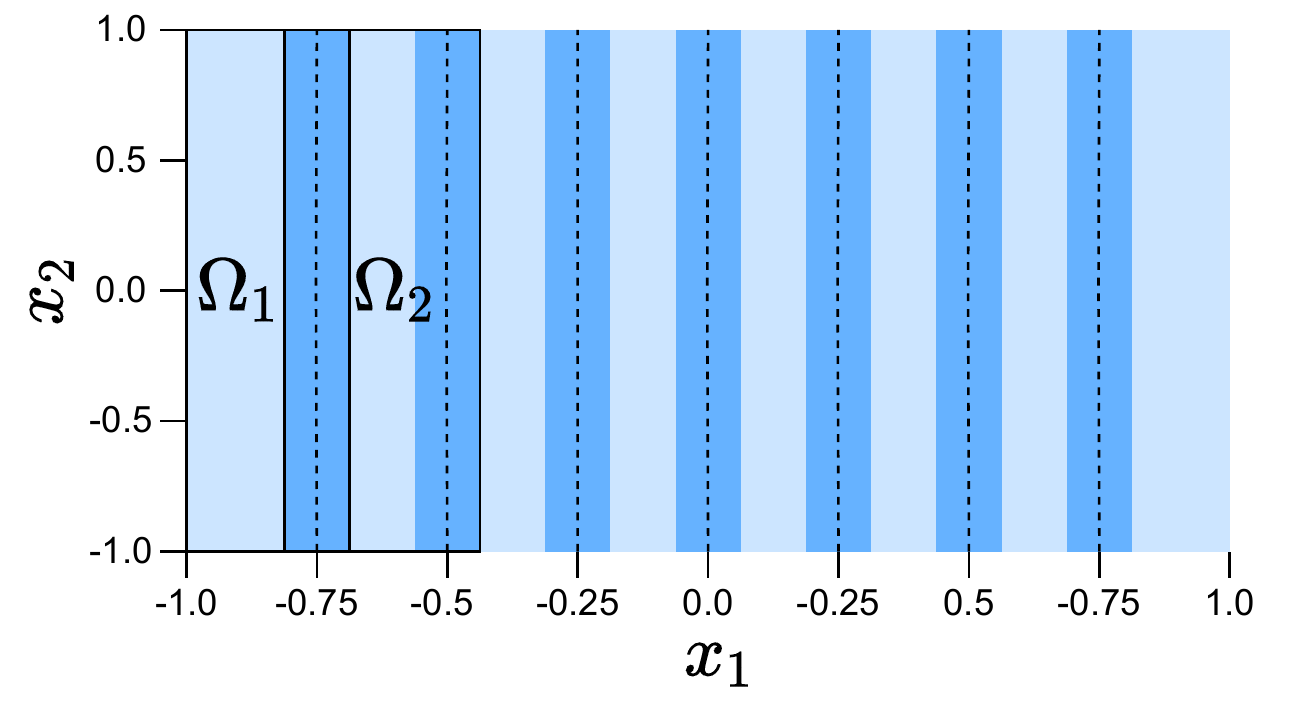}&
    \includegraphics[width=6cm]{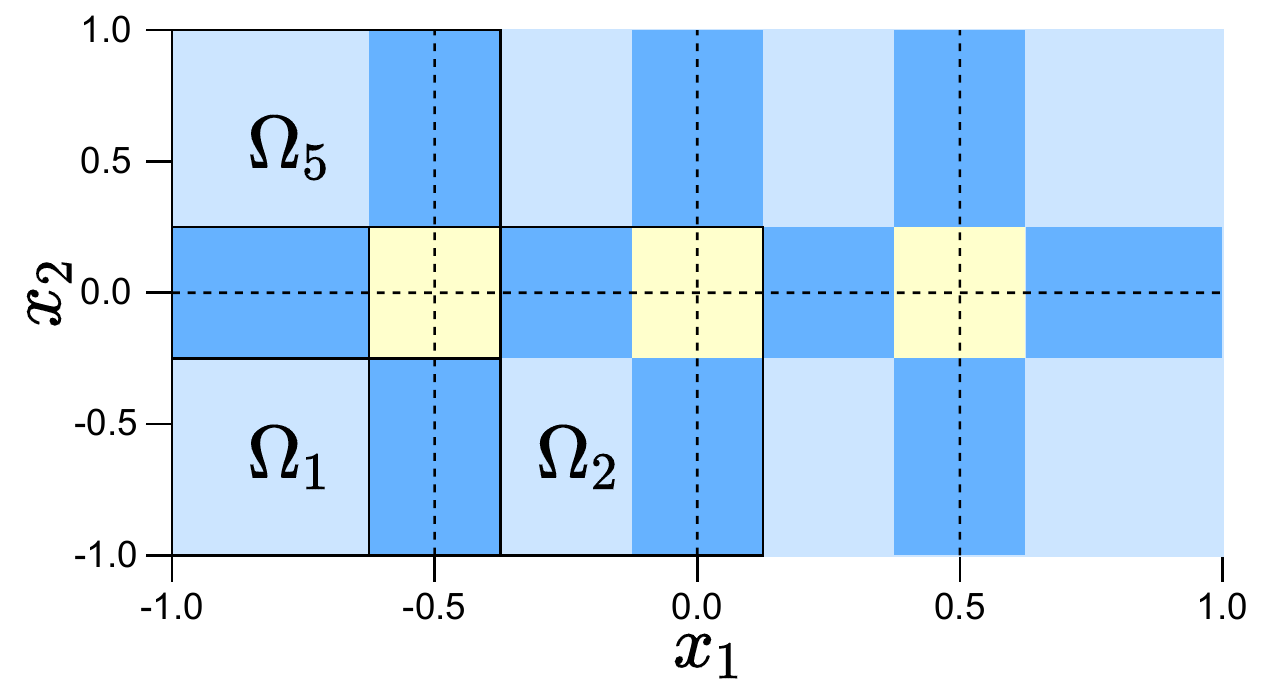}\\
    (a) $8\times 1$ subdomains & (b) $4\times 2$ subdomains\\
    \end{tabular}}
    \caption{
Two cases of domain decomposition with $N_d=8$ subdomains, two-dimensional Helmholtz equation. 
    }
    \label{fig_dd_8d}
\end{figure}

Two cases for decomposing the spatial domain  $\Omega$  into $N_d=8$ subdomains are considered, one where subdomains consist of 
vertical strips, the other using squares (see Figure \ref{fig_dd_8d}). 
For the first case (Figure \ref{fig_dd_8d}(a)), subdomains are given by
$\{\Omega_i\}_{i=1}^{N_d}=\{[x_{left}^{(i)},x_{right}^{(i)}]\times[-1,1]\}_{i=1}^{N_d}$, where $x_{left}^{(i)}$ and $x_{right}^{(i)}$ are defined in \eqref{x1} and  the length of the overlapping region is set to $w=0.125$.
For the second case (Figure \ref{fig_dd_8d}(b)), letting $N_d=N_{d1}\times N_{d2}$ with $N_{d1}=4$ and $N_{d2}=2$, each subdomain is defined as $\{\Omega_i\}_{i=1}^{N_d}=\{[x_{left}^{(i)},x_{right}^{(i)}]\times[x_{bottom}^{(i)},x_{top}^{(i)}]\}_{i=1}^{N_d}$ where  
\begin{equation}
    \begin{cases}
        x_{left}^{(i)}&=\max\left\{-1,-1+(i_1-1)\times\frac{2}{N_{d1}}-\frac{w_1}{2}\right\},\\
        x_{right}^{(i)}&=\min\left\{1,-1+i_1\times\frac{2}{N_{d1}}+\frac{w_1}{2}\right\},\\
        x_{bottom}^{(i)}&=\max\left\{-1,-1+(i_2-1)\times\frac{2}{N_{d2}}-\frac{w_2}{2}\right\},\\
        x_{top}^{(i)}&=\min\left\{1,-1+i_2\times\frac{2}{N_{d2}}+\frac{w_2}{2}\right\},
    \end{cases}
\end{equation}
where $i_1=i$ and  $i_2=1$ for $i=1,\ldots,4$,  $i_1=i-4$ and $i_2=2$ for $i=5,\ldots,8$,
the length of overlapping region along with $x_1$ is set to $w_1=0.25$ and that along with $x_2$ is set to $w_2=0.5$. 

In our F-D3M algorithm, MFFNets used for each subdomain contain two RFF subnetworks 
and the standard deviation $\sigma$ is set to $1$ and $10$ for the first and the second subnetworks respectively. Each RFF subnetwork has three hidden layers with thirty-two neurons. Similarly to the above two test problems, 
the setting of MFFNet in Global-MFFNet is the same as that for F-D3M and the network in Global-FCN is set to have three hidden layers with sixty-four neurons in each hidden layer. The weight coefficient in \eqref{em_loss_func} for Global-FCN and Global-MFFNet is set to $\lambda=100$. For all networks, the hyperbolic tangent function is employed as the activation function and the Xavier scheme is used for initialization \cite{glorot2010understanding}. All networks are trained by the Adam \cite{kingma2014adam} optimizer with default settings. An exponential learning rate with the initial learning rate of 0.01 and a decay rate of 0.9 every 1000 training epochs is employed, while for F-D3M, the learning rate is reinitialized for each domain decomposition iteration. As the local solutions in F-D3M gradually get close to the global solution as the domain decomposition iteration step
increases, we gradually increase the number of  training epochs---at the initial step, the number of epochs
is set to $2500$, and it increases by $500$ for each step until 
the number of epochs  reaches $1\times 10^5$.  For each training epoch, $|\mS_i|=625$ collocation points are generated with the uniform distribution with range $\Omega_i$ for $i=1,\ldots,N_d$. 
To train Global-FCN and Global-MFFNet, $10^5$ epochs are used and 
the number of collocation points for the interior domain and the boundary 
are set to $N_r=5000$ and $N_b=200\times4$ respectively.

\begin{figure}[!ht]
    \centerline{
    \begin{tabular}{c}
    \includegraphics[width=6cm]{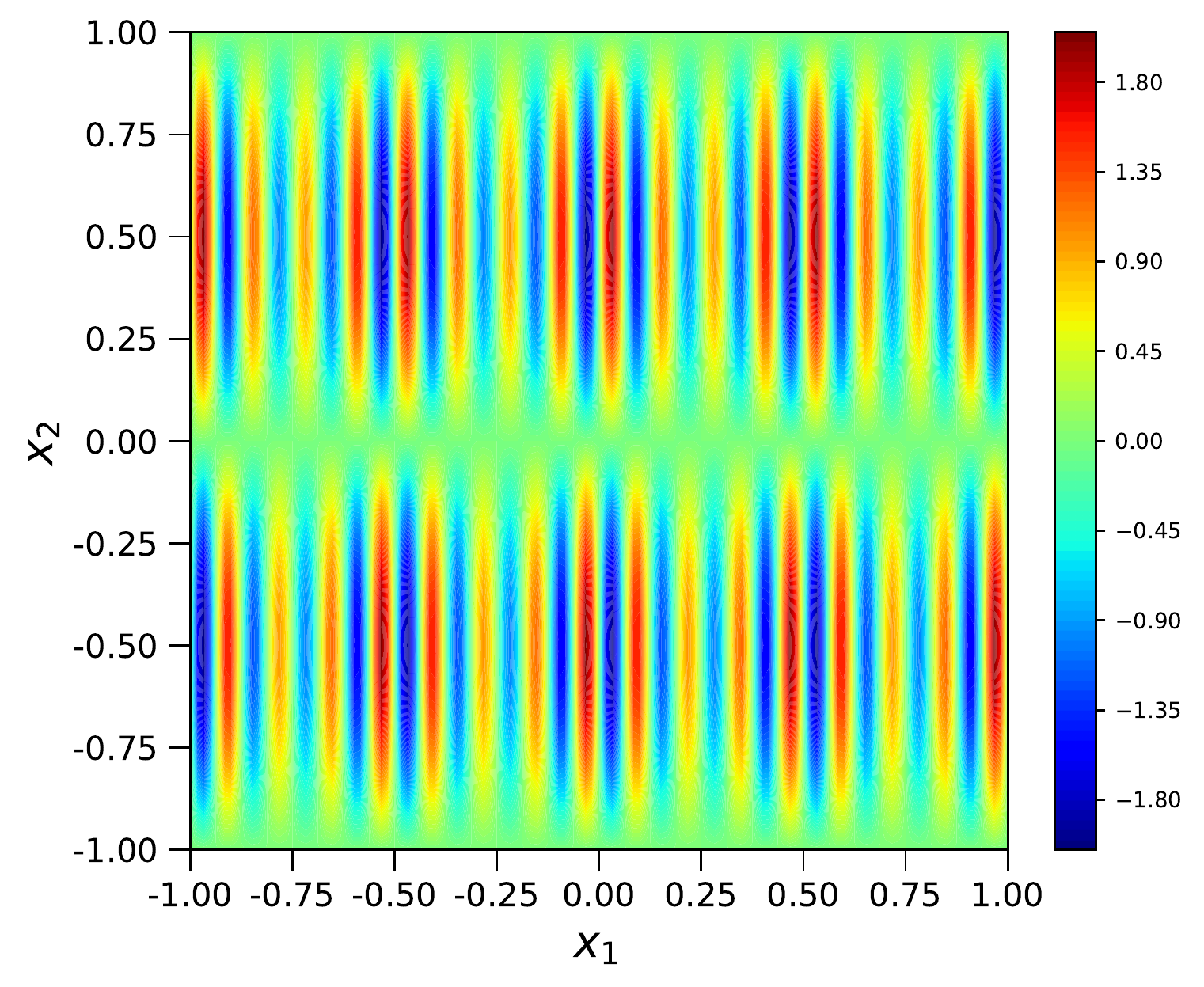}\\
    \end{tabular}}
    \caption{Exact solution, two-dimensional Helmholtz equation.}
    \label{fig_2d_Helmholtz_exact}
\end{figure}

\begin{figure}[!ht]
    \centerline{
    \begin{tabular}{ccc}
    \includegraphics[width=5cm]{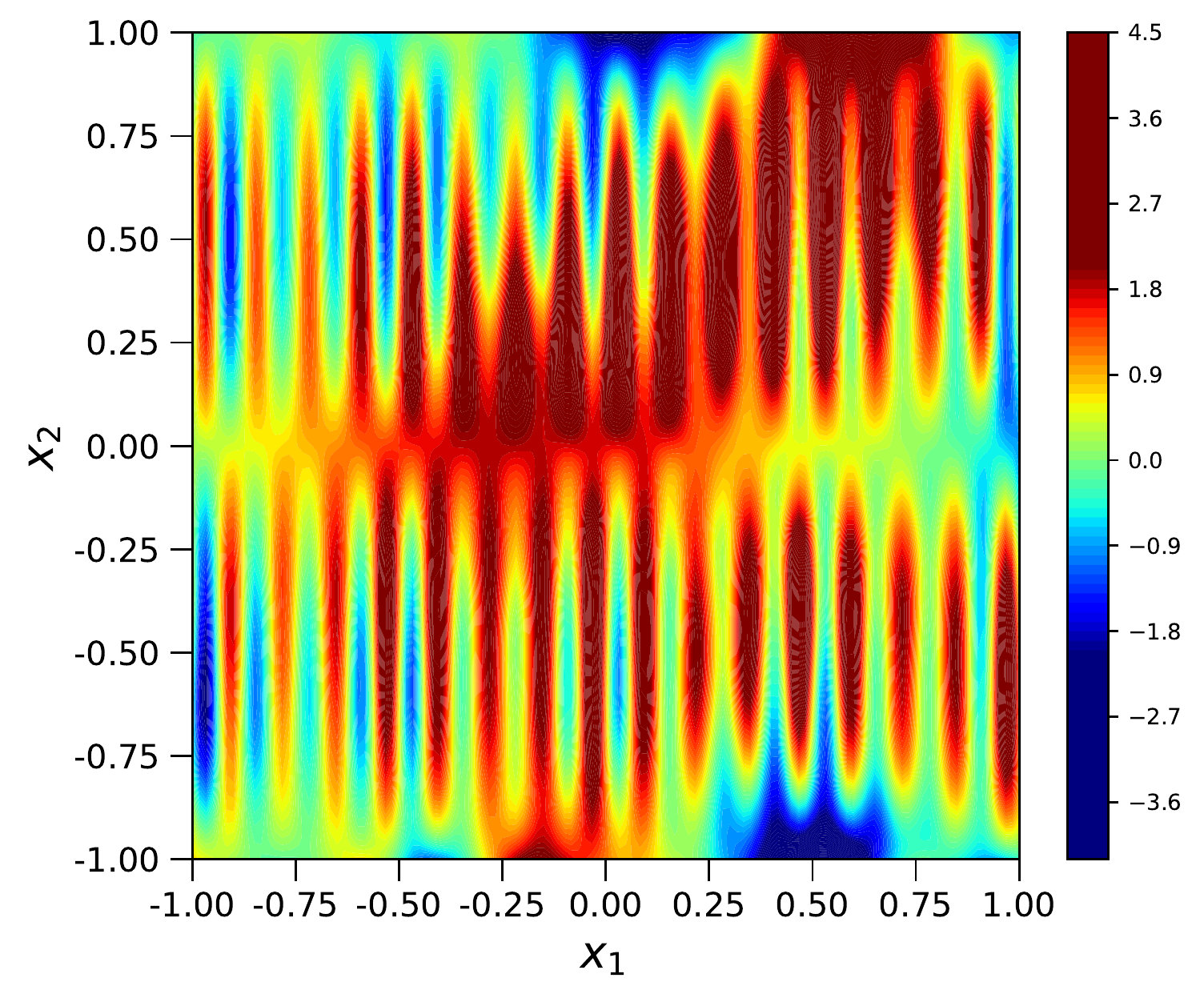}&\quad\quad\quad
    \includegraphics[width=5cm]{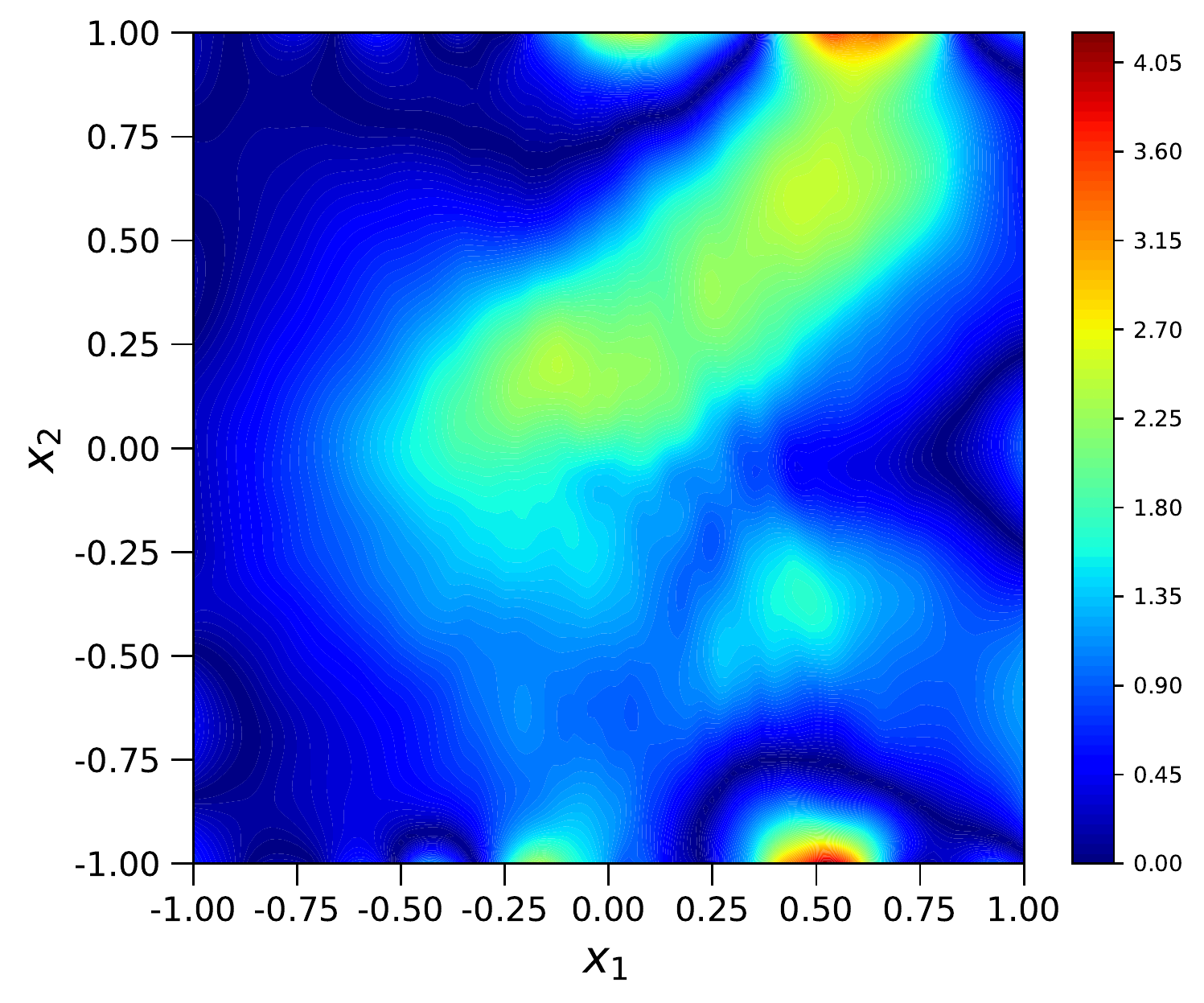}\\
    (a) Global-FCN appoximation &\quad\quad\quad (b) Global-FCN point-wise error\\
    \includegraphics[width=5cm]{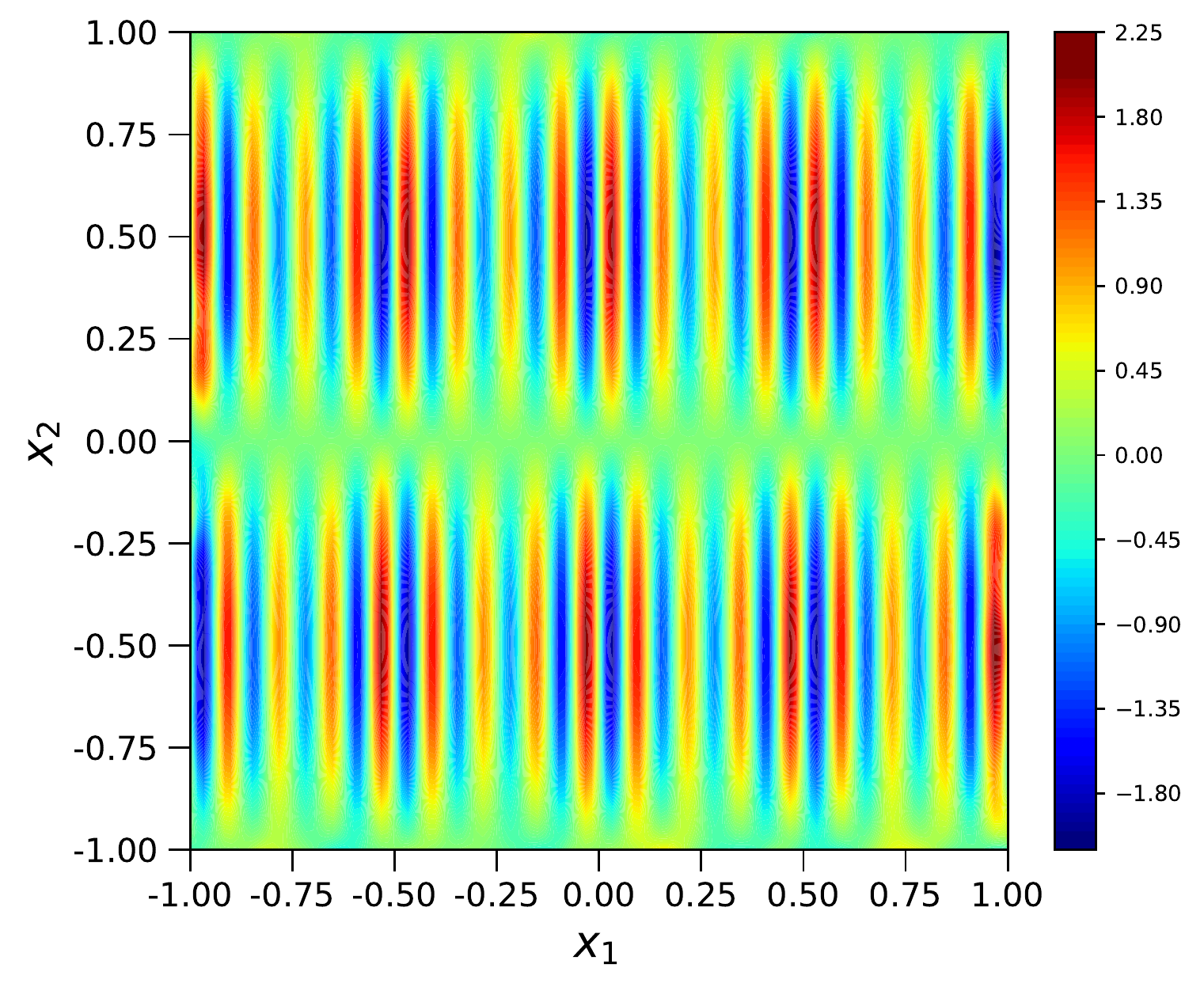}&\quad\quad\quad
    \includegraphics[width=5cm]{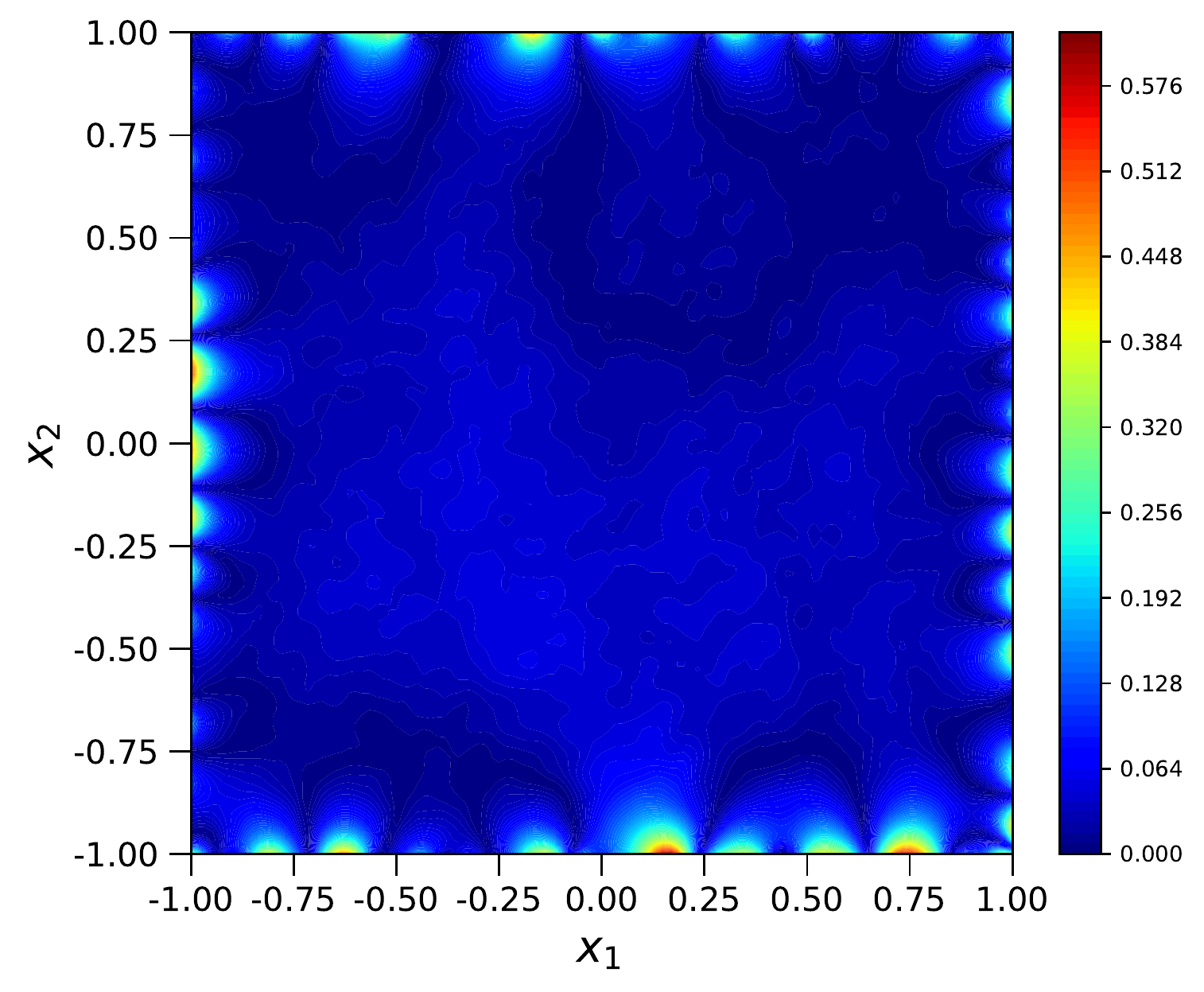}\\
    (c) Global-MFFNet appoximation &\quad\quad\quad (d) Global-MFFNet point-wise error\\
    \includegraphics[width=5cm]{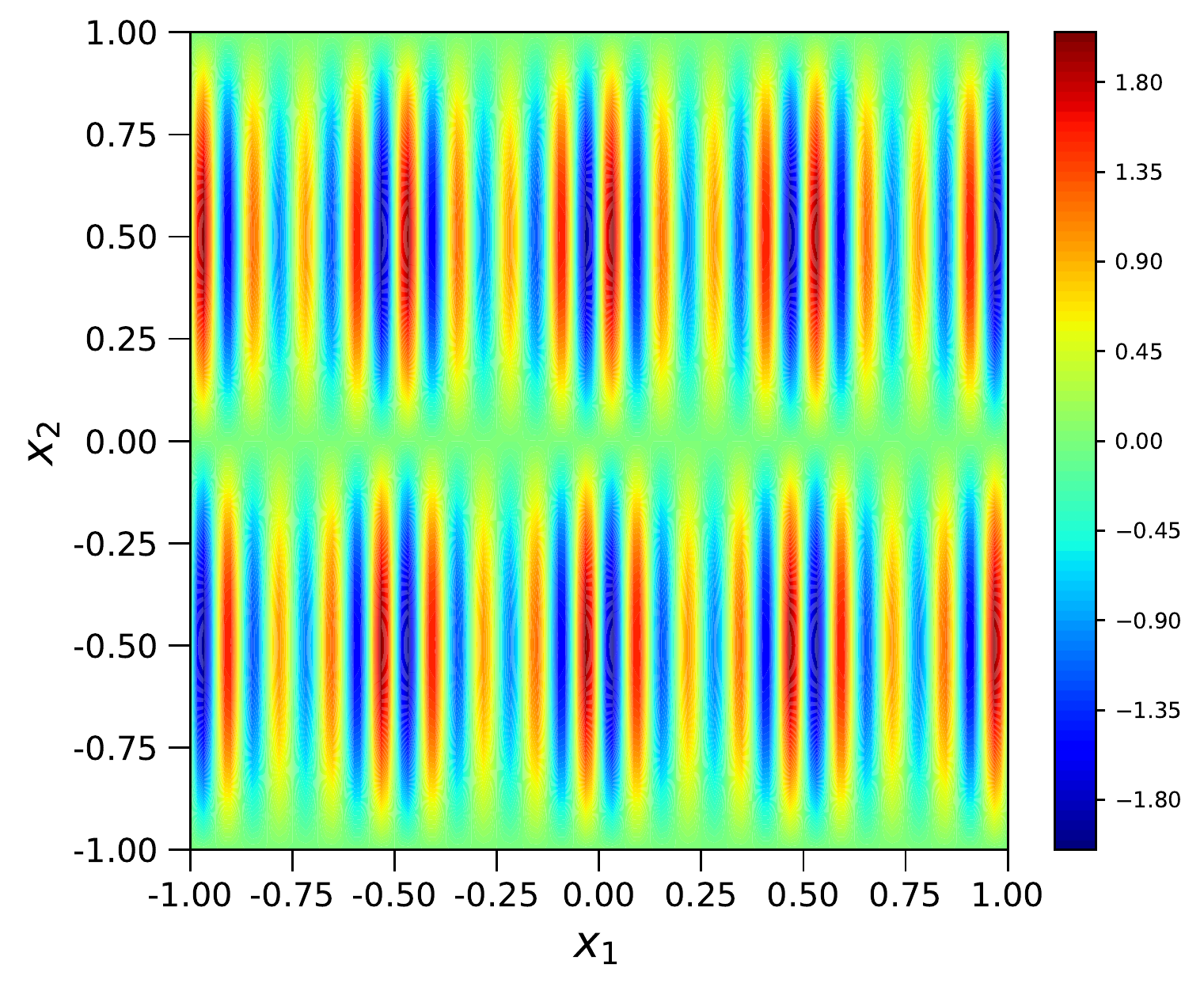}&\quad\quad\quad
    \includegraphics[width=5cm]{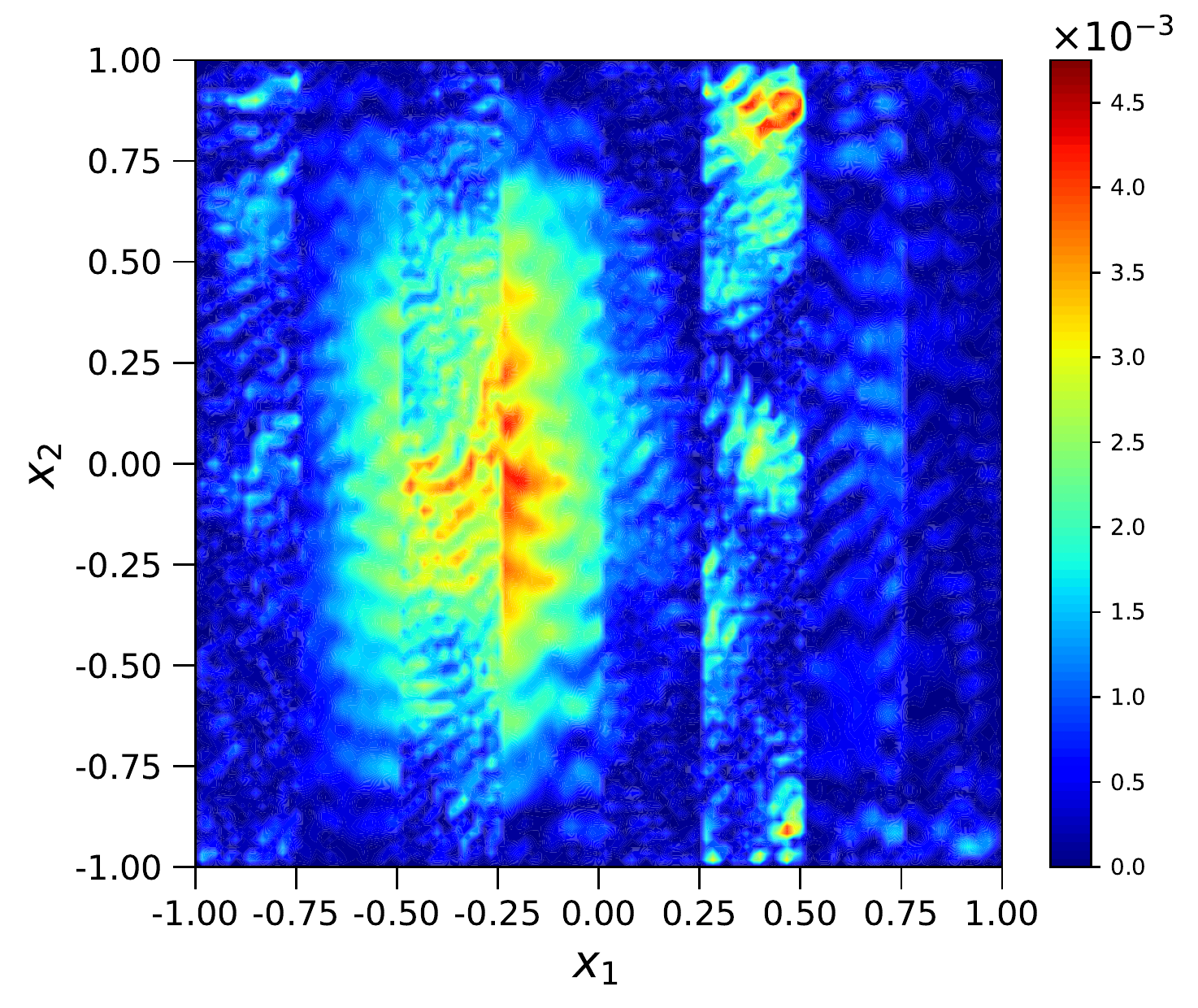}\\
    (e) F-D3M ($8\times 1$) appoximation &\quad\quad\quad (f) F-D3M ($8\times 1$) point-wise error\\
    \includegraphics[width=5cm]{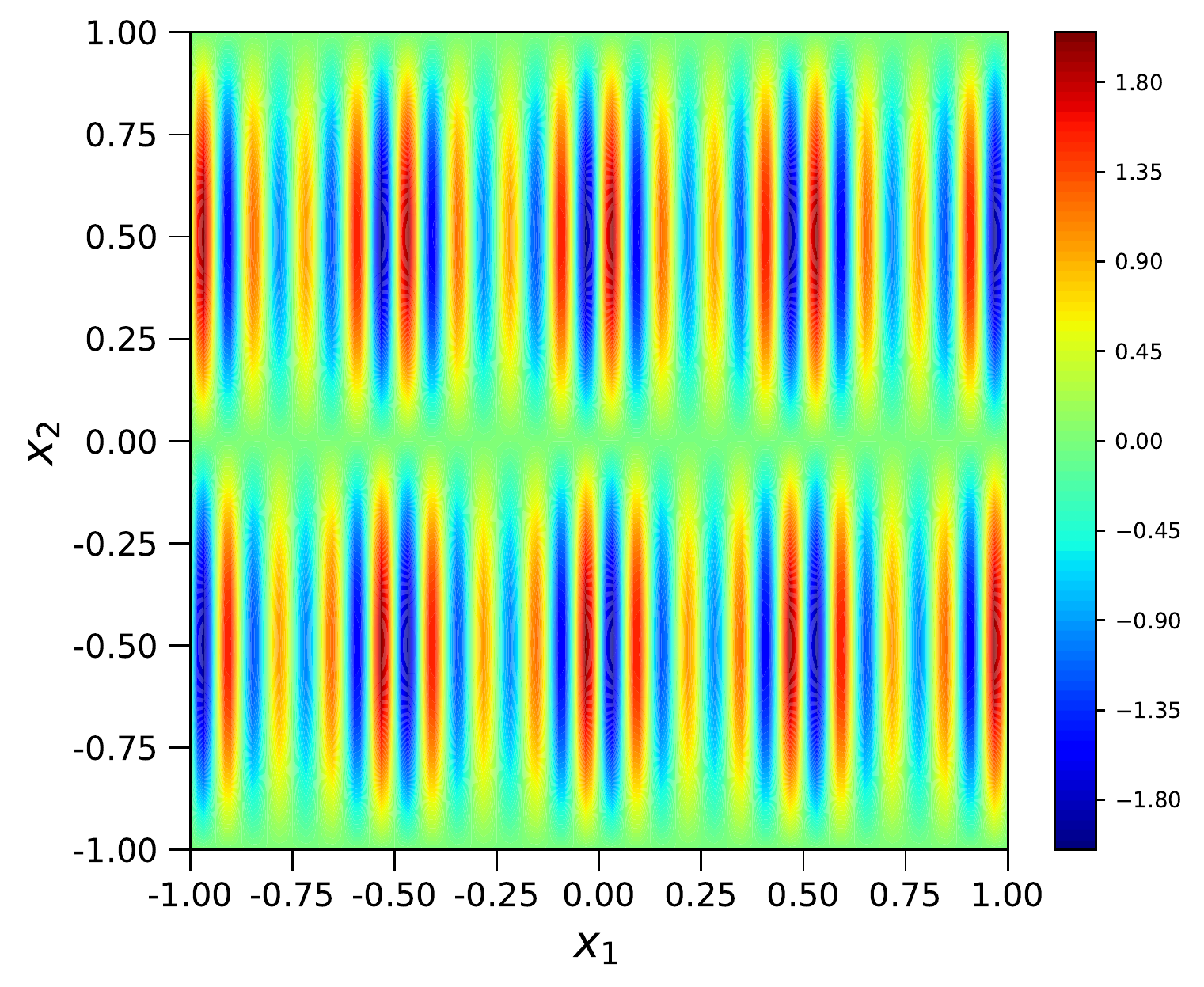}&\quad\quad\quad
    \includegraphics[width=5cm]{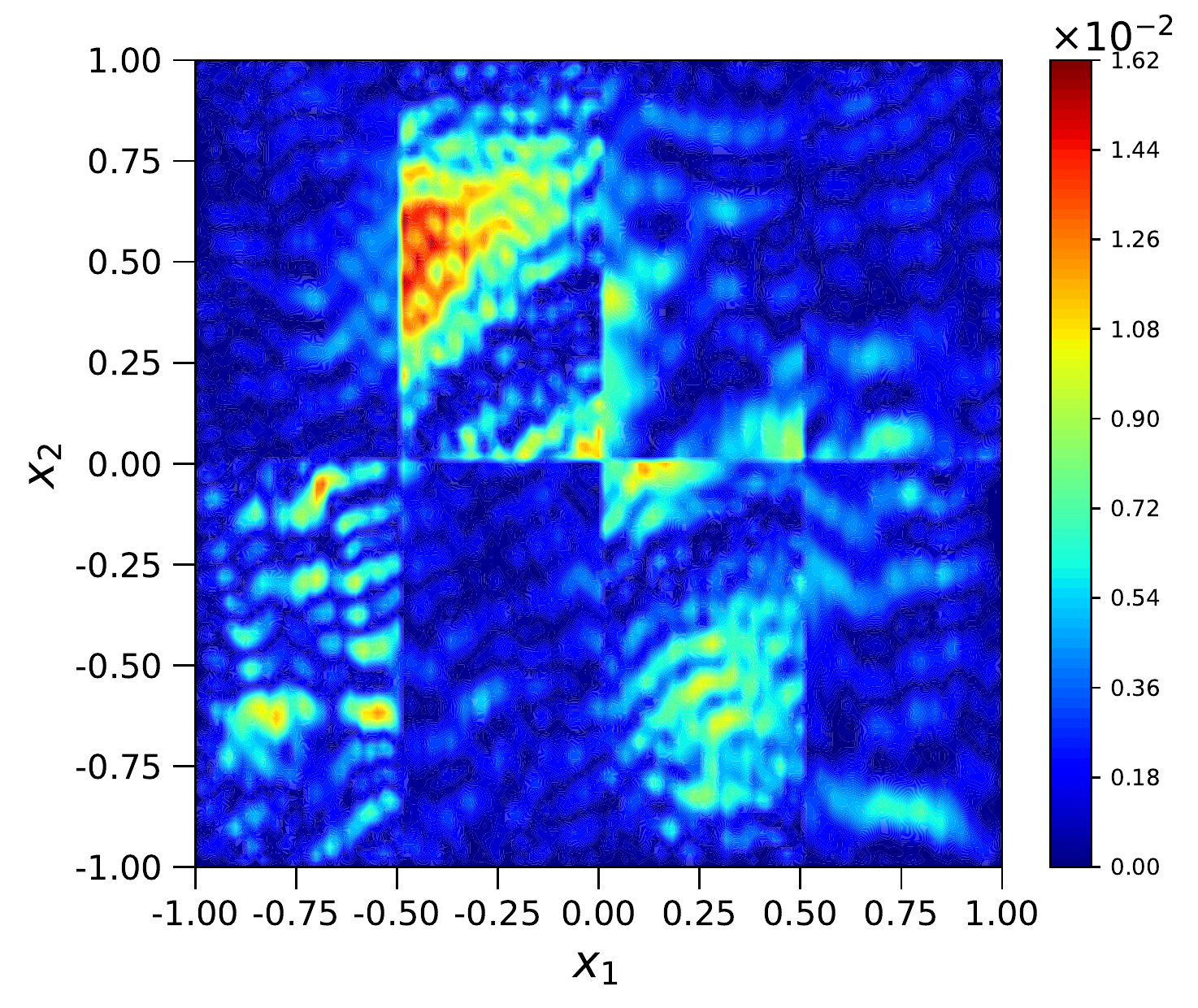}\\
    (g) F-D3M ($4\times 2$) appoximation &\quad\quad\quad (h) F-D3M ($4\times 2$) point-wise error\\
    \end{tabular}}
    \caption{\lr
    Approximation solutions and point-wise errors,  two-dimensional Helmholtz equation.
    }
    \label{fig_2d_Helmholtz_results}
\end{figure}

To assess the accuracy of approximation solutions,  
values of the exact solution and approximation solutions are computed 
for the uniform $121\times 121$ grids on $\Omega$. 
The absolute point-wise errors and the relative errors (see \eqref{relative_l2_error}) are computed, and all results are repeated three times for different random seeds to obtain average errors. 
Figure~\ref{fig_2d_Helmholtz_exact} shows the exact solution \eqref{test3-sol}, and Figure~\ref{fig_2d_Helmholtz_results} shows the approximation solutions and the corresponding absolute point-wise errors.  
It can be seen that the Global-FCN approximation is not close to the exact solution,
while there is no visual difference between approximations obtained by  F-D3M  (with both $8\times 1$ and $4 \times 2$ subdomains) and the exact solution. 
It also can be seen that the absolute point-wise errors of F-D3M are much smaller than the errors of 
Global-FCN and Global-MFFNet. 
The relative errors \eqref{relative_l2_error} of Global-FCN, Global-MFFNet, F-D3M with $8\times 1$ subdomains and F-D3M with $4\times 2$ subdomains are $6.63\mathrm{e}-01$, $9.56\mathrm{e}-02$, $1.54\mathrm{e}-03$ and $4.64\mathrm{e}-03$ respectively. 
It is clear that the approximation solutions of F-D3M are more accurate than the other two methods in accuracy---the relative error of F-D3M with $8\times1$ is two orders of magnitude smaller than that of Global-FCN and one order of magnitude smaller than that of  Global-MFFNet. 

\section{Conclusions} \label{sec:conclusions}
The divide and conquer principle is one of the fundamental concepts to solve PDE problems with high frequency 
modes. With a focus  on deep learning based methods, 
this paper proposes a Fourier feature based deep domain decomposition method (F-D3M). 
In F-D3M, difficulties caused by high frequency modes are curbed through  
decomposing global spatial domains into small overlapping subdomains, such that 
high frequency modes in the original problem can become relatively low frequency ones in local problems. 
In each subdomain, our multi Fourier feature network (MFFNet) provides effective approximations 
for local problem solutions. Especially, our results show that the training procedure of MFFNets is 
efficient with loss functions defined inside of subdomains. 
The overall efficiency of F-D3M is demonstrated with  numerical examples for diffusion and Helmholtz problems
with high frequency oscillations.
In the implementation F-D3M, the particular solution is updated during domain decomposition iterations, while MFFNets inherit the parameters from the last iteration directly. The mismatch between the particular solution and MFFNets can cause relative large errors during first few training epochs within each domain decomposition iteration. While proper initialization is crucial for training neural networks, optimal reinitialization strategies  for MFFNets in F-D3M still remain open challenges.  Developing efficient reinitialization strategies for MFFNets based on approximations obtained in  previous domain decomposition iterations will be the focus of our future work.

\bibliographystyle{elsarticle-num}
\bibliography{references}

\end{document}